\documentclass[10pt,english]{amsart}
\input{arxiv_preamble_v4}
\title{Residual supersingular Iwasawa theory and signed Iwasawa invariants}
\date{March $10^\text{th}$ 2021}
\author[Filippo A.~E.~Nuccio]{Filippo A. E. Nuccio Mortarino Majno di Capriglio}
\address{Univ Lyon, Université Jean Monnet Saint-Étienne \\ CNRS~UMR~5208, Institut Camille Jordan \\ F-42023 Saint-Étienne \\ France}
\email{filippo.nuccio@univ.st-etienne.fr}
\author{Ramdorai Sujatha}
\address{Department of Mathematics \\ University of British Columbia \\ Vancouver, BC~V6T~1Z2 \\ Canada}
\email{sujatha@math.ubc.ca}
\begin{document}
\begin{abstract}
For an odd prime $p$ and a supersingular elliptic curve over a number field, this article introduces a multi-signed residual Selmer group, under certain hypotheses on the base field. This group depends purely on the residual representation at $p$, yet captures information about the Iwasawa theoretic invariants of the signed $p^\infty$-Selmer group that arise in supersingular Iwasawa theory. Working in this residual setting provides a natural framework for studying congruences modulo $p$ in Iwasawa theory.
\end{abstract}
\maketitle
\section{Introduction}
Iwasawa theory of Galois representations, especially those arising from elliptic curves and modular forms, affords deep insights into the arithmetic  of these objects over number fields. The  Iwasawa theoretic invariants, especially the~$\iwmu[\empty]$ and~$\iwlambda[\empty]$ invariants, play a central role in this study.  Iwasawa theory for ordinary elliptic curves, and more generally for ordinary Galois representations, was initiated by Mazur in~\cite{Maz72} and Greenberg in~\cite{Gre89}. The corresponding theory for supersingular elliptic curves is subtler and was already begun by Perrin-Riou in~\cite{Per90}. In the last couple of decades, supersingular Iwasawa theory has gained considerable momentum (see~\cite{Kob03,Pol03,IovPol06,LeiLoeZer10,Spr12,Kim13,Kim18,KitOts18} and references therein). 

Greenberg and Vatsal \cite{GreVat00} investigated the behaviour of Iwasawa invariants for ordinary elliptic curves whose residual representations are congruent. The objects of study are the dual $p^\infty$-Selmer groups of the elliptic curves over the cyclotomic $\ZZ_p$-extension of the base field, which is assumed to be a number field.
Specifically, let $p$ be an odd prime and  $\EC_i$, $i=1,2$ be two elliptic curves over $\QQ$ with good ordinary reduction at $p$. Greenberg and Vatsal prove that the vanishing of the $\iwmu[\empty]$-invariant for the dual $p^\infty$-Selmer group of one of the curves implies the vanishing for the other. Their study makes crucial use of a \emph{non-primitive} dual Selmer group, which has the same $\iwmu[\empty]$-invariant as the dual $p^\infty$-Selmer group.
When the $\iwmu[\empty]$-invariants vanish, they also prove the equality of the $\iwlambda[\empty]$-invariants for the  non-primitive dual $p^\infty$-Selmer groups for $\EC_1$ and $\EC_2$. However, they provide examples showing that the $\iwlambda[\empty]$-invariants for the dual $p^\infty$-Selmer groups do not coincide. These results have been extended to the representations coming from higher weight modular forms by Emerton, Pollack and Weston in~\cite{EmePolWes06}, and to more general base fields and $\ZZ_p$-extensions by, among others, Hachimori in~\cite{Hac11} and by Kidwell in~\cite{Kid18}. A  crucial input in the study of $p^\infty$-Selmer groups in the ordinary case is a deep result of Kato (see~\cite{Kat04}) which implies that the dual $p^\infty$-Selmer groups (and their non-primitive counterparts) are torsion modules over the Iwasawa algebra.

When $\EC/\QQ$ is an elliptic curve having good, supersingular reduction at $p$, the dual $p^\infty$-Selmer group is no longer torsion over the Iwasawa algebra. Kobayashi defined the \emph{signed} $p^\infty$-Selmer group in~\cite{Kob03}, making use of special subgroups of the local Mordell--Weil groups along the cyclotomic tower which were already considered by Perrin-Riou. These signed $p^\infty$-Selmer groups are torsion over the Iwasawa algebra and display properties that are strikingly similar to those of the $p^\infty$-Selmer group in the ordinary case and come equipped with \emph{signed} Iwasawa invariants $\iwlambda[\pm],\iwmu[\pm]$. Results analogous to those of Greenberg--Vatsal for these signed invariants were proved by Kim in~\cite{Kim09}, again making use of the  non-primitive Selmer groups. The study of signed Selmer groups for higher weight modular forms has been initiated by Lei, Loeffler and Zerbes in~\cite{LeiLoeZer10} through the theory of Wach modules, and extensions of Greenberg--Vatsal results in this setting can be found in~\cite{HatLei19} by Hatley and Lei. The definition of the signed Selmer groups has been extended to a broader class of number fields in~\cite{IovPol06,Kim13,KitOts18}. In this article, we will mainly refer to Kitajima--Otsuki's paper. Our work sheds more light on the behaviour of the Iwasawa invariants for the dual signed $p^\infty$-Selmer groups of elliptic curves in the supersingular case. The results proved here are more general than
those in~\cite{Kim09}, largely  because we set up  a  framework for multi-signed Selmer groups.

The novelty in our approach is that we systematically work with the residual representation of a supersingular elliptic curve defined over a number field $\basefield$ satisfying certain conditions (see~Section~\ref{section:preliminaries}, in particular hypothesis~\ref{hyp_1} therein), instead of working with $\tors{\EC}{p^\infty}$. In particular, we introduce a new Selmer group, attached to the Galois representation $\tors{\EC}{p}$ of $p$-torsion points of $\EC$, which we call multi-signed residual Selmer group. It depends only on the isomorphism class of the residual Galois representation $\tors{\EC}{p}$, yet captures the full Iwasawa-theoretic information about the $\iwmu[\pm]$-~and the $\iwlambda[\pm]$-invariants of the usual signed $p^\infty$-Selmer group. The multisigned residual Selmer group contains the residual analogue of the fine Selmer group as defined  by Coates and the second author in ~\cite{CoaSuj05}.
.In~\cite{CoaSuj05}. \emph{Ibid.} the authors postulate a conjecture, referred to as~\ref{conj:A}, which asserts that the Iwasawa $\iwmu[\empty]$-invariant of the dual fine $p^\infty$-Selmer group over $\basecyc$ vanishes. It is pertinent to remark here that~\ref{conj:A} depends only on the residual Galois representation~(see~\cite{Gre11} and~\cite{Suj10}) and its formulation is independent of the reduction type at $p$ of the elliptic curve. Working directly with the multi-signed residual Selmer group provides a conceptual framework to explore the comparison of Iwasawa invariants, when the residual representations are isomorphic. It also potentially provides the right context for explaining a plethora of congruences in arithmetic, such as the congruences between complex and $p$-adic $L$-values which occur when the residual representations are isomorphic. We hope to return to this subject of framing a \emph{residual Iwasawa theory} in our future works.

The main results of this paper are Theorem~\ref{thm:main_result} and Theorem~\ref{thm:mu=mu}. Under certain hypothesis~\ref{hyp_1} and~\ref{hyp_2}, and assuming~\ref{conj:A}, Theorem~\ref{thm:main_result} provides a criterion for the $\iwmu[\empty]$-invariant of the multi-signed $p^\infty$-Selmer group to vanish, purely in terms of the multi-signed residual Selmer group. We refer to the main body of the paper for its statement, because it involves some morphism whose definition is too technical for this introduction.

As an application of Theorem~\ref{thm:main_result}, our second theorem provides a criterion for the $\iwmu[\empty]$-invariant of the signed $p^\infty$-Selmer group to vanish, purely in terms of the multi-signed residual Selmer group. In the following, denote by $\ordssprimeset$ the number of primes in $\basefield$ of supersingular reduction for an elliptic curve $\EC$ and let $\multisign\in\{+,-\}^{\ordssprimeset}$ be a vector of signs, in other words $\multisign$ is a $\ordssprimeset$-tuple with entries either $+$ or $-$. 
Denote by $\KobdualSel*{\tors{(\EC)}{p^\infty}}{\basecyc}$ the dual multi-signed $p^\infty$- Selmer groups, as defined in~Definition~\ref{def:pm-Sel} (see also \cite[Definition~2.1]{KitOts18}):
\begin{introtheorem}{\ref{thm:mu=mu}}
Let $\EC_1,\EC_2$ be two elliptic curves defined over $\basefield$, satisfying hypothesis~\ref{hyp_1} and such that the residual Galois representations $\tors{(\EC_1)}{p}$ and $\tors{(\EC_2)}{p}$ are isomorphic. Then, the sets $\ssprimeset[1]$ and $\ssprimeset[2]$ of primes of supersingular reduction for $\EC_1$ and $\EC_2$ coincide. Given a vector $\multisign\in\{+,-\}^{\ordssprimeset}$, assume that both curves satisfy~\ref{hyp_2} and let $\iwmu_{\EC_j}$ be the Iwasawa $\iwmu[\empty]$-invariants of $\KobdualSel*{\tors{(\EC_j)}{p^\infty}}{\basecyc}$, for $j=1,2$. Then 
\begin{equation*}
\iwmu_{\EC_1}=0\Longleftrightarrow \iwmu_{\EC_2}=0.
\end{equation*}
\end{introtheorem}
In the ordinary case,  using the non-primitive Selmer groups, Greenberg-Vatsal (see~\cite[Theorem~1.5]{GreVat00}),  compare the $\iwlambda[\empty]$-invariants. When the $\iwmu[\empty]$-invariants vanish, our approach also  enables 
such a comparison, without the  use of the imprimitive Selmer groups. These results should be compared with those by Kim (see~\cite[Corollary~2.13]{Kim09}) and Hatley--Lei (see~\cite[Theorem~4.6]{HatLei19}). We need the  additional assumption that the Pontryagin dual of the usual Selmer group (and hence, also of the multi-signed one by Corollary~\ref{cor:no_finite_in_signedSel}) contains no non-zero finite $\IwLambda$-submodule, and this is known to be true in many cases.
\begin{introtheorem}{\ref{thm:lambda_inv}}
Let $\EC_1,\EC_2$ be two elliptic curves defined over $\basefield$, satisfying the hypotheses of Theorem~\ref{thm:mu=mu}. Assume that their Iwasawa $\iwmu$-invariants vanish and suppose further that the Pontryagin duals $\dualSel*{\tors{(\EC_j)}{p^\infty}}{\basecyc}$ of the usual Selmer groups do not have any non-zero finite $\IwLambda$-submodule. Then,
\begin{equation*}
\iwlambda_{\EC_j}=\seldim{\empty}+\defect{\EC_j}
\end{equation*}
where $\defect{\EC_j}$ is as in Definition~\ref{def:defect} and $\seldim{\empty}:=\seldim{\EC_1}=\seldim{\EC_2}$ is as in~\eqref{eq:def_rho}.
\end{introtheorem}

The first term $\seldim{\empty}$ in the above statement depends only on the residual representation and is independent of the curve. The second term $\defect{\EC_j}$ depends on the structure of the local $p$-torsion of the elliptic curve over the first layer of the cyclotomic tower at the set of primes of bad reduction, together with the primes above $p$ with ordinary reduction: in particular, it is independent of the vector~$\multisign$. The reader is referred to the main body of the paper for the precise definitions of these numerical invariants.

Our methods also show that the difference $\iwlambda[\multisign]-\iwlambda[-\multisign]$ depends only on the residual representation, a fact which was already observed by Kim in~\cite[Remark~3.3]{Kim09} for $\basefield=\QQ$. Propositions~\ref{prop:universal_sequence} and ~\ref{prop:raynaud} are the key technical tools needed to show that the definition of the multi-signed residual Selmer group depends only on the residual representation. They compare the reduction type at places above $p$ of two residually isomorphic elliptic curves at primes above $p$, and are of independent interest. While the first result relies on Honda--Tate theory and is a typical feature of supersingular reduction,  the second relies upon a result by Raynaud on finite flat group schemes killed by $p$. We remark that some of the techniques used here were first developed in the setting of purely ordinary reduction in~\cite{LimSuj18a}. The extension to the generality considered in this paper was not obvious then.

In ongoing works, we extend these results in the following different directions. First, to higher weight modular forms over the cyclotomic extension, second to multiple $\ZZ_p$-extensions, and finally to the multi-signed Selmer groups as well as non-commutative $p$-adic Lie extensions.

The paper consists of five sections, including this introductory section. In Section~\ref{section:preliminaries}, we introduce notation and some preliminaries about the local structure of elliptic curves with supersingular reduction. In Section~\ref{section:signed}, we recall the main properties of multi-signed Kummer maps and Selmer groups, mainly building upon~\cite{KitOts18}, and we introduce the multi-signed residual Selmer group. In Section~\ref{section:iwasawa}, we study the Iwasawa theory of this group and state our main results. The final section presents some numerical examples that illustrate our results.

\section*{Acknowledgment} We would like to thank Laurent Berger for inviting the second author to ENS Lyon in June~$2018$ which led to initiating this research. This work was continued during the visit of the first author to the Pacific Institute of Mathematical Sciences (PIMS), Vancouver benefiting of a CNRS-PIMS exchange. F.~A.~E.~N. gratefully acknowledges the support and hospitality of PIMS and of the University of British Columbia, as well as the \emph{acceuil en délégation à l'UMI-3069 du CNRS}. S.~R. gratefully acknowledges support from NSERC Discovery grant 201903987. We would also like to thank Matthieu Romagny for helpful correspondence. 
\section{Preliminaries}\label{section:preliminaries}
In this paper $\basefield$ denotes a fixed number field of absolute degree $[\basefield:\QQ]=\degbase$, and $\EC/\basefield$ is an elliptic curve
defined over $\basefield$. Throughout, $p$ will denote an odd prime $\geq 3$, $\pprimeset$ denotes the set of primes above $p$ in $\basefield$ and $\tate{\EC}$ will denote the Tate module of $\EC$. The following hypothesis, which will be referred to as~\ref{hyp_1}, is assumed throughout (\cfr~\cite[Theorem~1.3~(i)--(v)]{KitOts18}):
\begin{equation}\tag*{\textup{\textbf{Hyp~1}}}\label{hyp_1}
\begin{minipage}{0.9\textwidth}
\begin{enumerate}[label=\textup{(}\roman*\textup{)}]
\item The curve $\EC/\basefield$ has good reduction at all primes in $\pprimeset$;
\end{enumerate}
Denote by $\ssprimeset\subseteq \pprimeset$ the set of primes above $p$ where $\EC$ has supersingular reduction.
\begin{enumerate}[resume,label=\textup{(}\roman*\textup{)}]
\item $\ssprimeset$ is non-empty;
\item all primes $\ssprimeid{1},\hdots,\ssprimeid{\ordssprimeset}$ in $\ssprimeset$ split completely in $\basefield/\QQ$, so $\genlocal{\ssprimeid{i}}\cong\QQ_p$ for all $1\leq i\leq \ordssprimeset$;
\item $1+p-\lvert{\modEC(\finitefield[\ssprimeid{i}])}\rvert=0$, where $\modEC$ is the reduction of $\EC$ modulo any of the prime ideals $\ssprimeid{1},\hdots,\ssprimeid{\ordssprimeset}$ and $\finitefield[\ssprimeid{i}]$ denotes the residue field of $\ssprimeid{i}$.
\item The ramification index $e(\ordprimeid{})$ in the extension $\basefield/\QQ$ of every prime $\ordprimeid{}\in\pprimeset$ where $\EC$ has good, ordinary reduction, is at most $p-1$.
\end{enumerate}
\end{minipage}\end{equation}
\begin{remark} Kitajima and Otsuki work in a slightly greater generality, allowing the supersingular primes to be simply unramified in $\basefield/\QQ$, provided the curve is defined over a subfield where they split completely. Points (i)--(iv) in~\ref{hyp_1} ensure that the signed Selmer groups are defined (see Section~\ref{section:signed}). Point~(v) in~\ref{hyp_1} ensures that the multi-signed residual Selmer group depends only on the residual Galois representation $\tors{\EC}{p}$ (\cfr Proposition~\ref{prop:universal_sequence} and Proposition~\ref{prop:raynaud}).
\end{remark}
Let us set the notation that will be used in the paper.
\begin{notation}\label{not:S_and_primes}
The set $\pprimeset$ of primes above $p$ is the disjoint union $\pprimeset=\ssprimeset\coprod \ordprimeset$, where $\ordprimeset=\{\ordprimeid{1},\hdots,\ordprimeid{\ordordprimeset}\}$ is the (possibly empty) set of primes where $\EC$ has good, ordinary reduction. Consider the cyclotomic $\ZZ_p$-extension $\basecyc/\basefield$ of $\basefield$, with intermediate layers $\basefield[n]$, for $n\geq 0$, so that $\Gal{\basefield[n]}{\basefield}\cong\ZZ/p^n\ZZ$ and $\basefield[0]=\basefield$.  
By~\ref{hyp_1}, all primes $\ssprimeid{i}\in\ssprimeset$ split in $\basefield$, so they are all totally ramified in $\basecyc/\basefield$. Let $\ssprimeid[n]{i}$ denote the prime ideal of $\basefield[n]$ above $\ssprimeid{i}$, for $1\leq i\leq \ordssprimeset$ and let $\localfield[i]{n}$ be the localisation of $\basefield[n]$ at $\ssprimeid[n]{i}$. For ease of notation, we often suppress the index $i$ since these fields, for fixed $n$, are all isomorphic to the $n$th layer of the cyclotomic extension of $\localfield{0}\cong\QQ_p$. In particular, $\localfield[i]{\cyc}\cong(\Qcyc)_{\ssprimeid[\cyc]{i}}$. We also need to consider the fields obtained by adjoining $p$-power order roots of unity to $\localfield[i]{0}=\localfield{0}$. Set $\fullcyc{n}=\localfield{0}(\zeta_{p^{n+1}})$, for $n\geq -1$, where $\zeta_{p^{n+1}}$ is a primitive $p^{n+1}$th root of unity. For all $n\geq -1$, we let $\maxid{n}$ be the maximal ideal of $\fullcyc{n}$. In particular, $\fullcyc{-1}\cong\QQ_p$ and $\maxid{-1}\cong p\ZZ_p$. The Galois group $\Gal{\fullcyc{0}}{\localfield{0}}$ is denoted by $\galprime$. It is isomorphic to $\Gal{\fullcyc{n}}{\localfield{n}}$ for all $n\geq 0$, and we tacitly identify these groups throughout.

Let $\primeset=\pprimeset\coprod\badprimeset$ where $\badprimeset=\{\badprimeid{1},\hdots,\badprimeid{\ordbadprimeset}\}$ is the finite set of primes of bad reduction for $\EC/\basefield$. The maximal extension of $\basefield$ unramified outside of $\primeset$ will be denoted $\baseS$. We usually write $\genprimeid$ or $\GenPrimeId$ to denote generic primes above $\primeset$ in an extension of $\basefield$. Given an extension $\basefield'/\basefield$, we sometimes abuse notation and again denote by $\primeset$ the primes of $\basefield'$ that lie above primes in $\primeset$. When we need to specify the field, we write $\genprimeset[\basefield']$ for $\ast\in\{\emptyset,\mathrm{ord},\mathrm{ss},\mathrm{bad}\}$ to denote the sets of primes of $\basefield'$ above primes in $\genprimeset$.

Given any field $\field\in\{\basefield[n],\genlocal[n]{\genprimeid},\fullcyc{n}\}$ (for some $0\leq n<\infty$ and possibly some prime $\genprimeid\in \primeset$ of $\basefield[n]$), its ring of integers will be denoted by $\rint{\field}$; when $\field$ is a local field, we further denote its residue field by $\finitefield[\genprimeid]$.
For $\field$ as above, write $\fieldS=\baseS$ if $\field=\basefield[n]$, $\fieldS=\overline{\fullcyc{n}}=\overline{\QQ_p}$ if $\field=\fullcyc{n}$ and $\fieldS=\baseSloc{\GenPrimeId}=\overline{\genlocal[n]{\genprimeid}}$, for some extension $\GenPrimeId\mid \genprimeid$, when $\field=\genlocal[n]{\genprimeid}$. The corresponding Galois groups $\Gal{\fieldS}{\field}$ are denoted, respectively, by $\GS[n]$, $\Glocal[\fullcyc]{n}$ and $\Glocal[{\genlocal[n]}]{\genprimeid}$; 
in case $\genprimeid=\ssprimeid{i}$, this will be denoted $\Glocal{i}$. When $\field\in\{\genlocal[n]{\genprimeid},\fullcyc{n}\}$, and $M$ is any Galois module, we usually write $\HH{i}{\field}{M}$ to denote the cohomology group $\HH{i}{\Gal{\fieldS}{\field}}{M}$.
\end{notation}
The Galois module of $p^t$-torsion points of $\EC$ is denoted $\tors{\EC}{p^t}$, and more generally $\tors{M}{p^t}$ will denote the submodule of a Galois module $M$ consisting of the $p^t$-torsion elements in $M$. By a slight abuse of notation, $\modEC/\finitefield$ is the reduction of $\EC$ modulo any of the prime ideals $\ssprimeid{i}$, all the reductions being isomorphic. Similarly, $\formalEC$ is the formal group of $\EC$ over $\ZZ_p=\rint{\localfield{0}}$. As discussed in \cite[Corollary~8.5]{Kob03} and \cite[\S~3.1]{KitOts18}, there is a $\ZZ_p$-isomorphism $\formalEC\cong\Fss$ between the formal group of $\EC/\ZZ_p$ and the supersingular formal group $\Fss$ whose logarithm is of Honda type $t^2+p$ (the group $\Fss$ is denoted by $\formalKO$ in~\cite[\S 3.1]{KitOts18}: observe that in our setting the automorphism $\varphi$ in \loccit is trivial).
\section{Plus and minus decomposition}\label{section:signed}
\subsection{The signed Kummer maps}
The aim of this section is to gather some results about the plus/minus decomposition of \emph{local} Mordell--Weil and formal groups, mainly taken from~\cite{KitOts18}, which in turn relies on~\cite{Kob03}. Most of the results mentioned below are either well-known or easy adaptations to the finite Galois module of $p$-torsion points, of arguments which are normally stated for the divisible module of $p^\infty$-torsion points.

We start with a general remark about vanishing of global torsion points for $\EC$ along the cyclotomic extension.

\begin{proposition}\label{prop:no_p_torsion} For every $n\geq 0$, the torsion subgroup $\tors{\formalEC(\maxid{n})}{p}$ is trivial. In particular,
\[
\tors{\EC(\fullcyc{n})}{p}=\tors{\EC(\localfield{n})}{p}=\tors{\EC(\basefield[n])}{p}=\{0\}\qquad\text{ for all }n\geq 0.
\]
\end{proposition}
\begin{proof} See {\cite[Proposition~3.1]{KitOts18}}.
\end{proof}

Following the pivotal works \cite{Per90} and \cite{Kob03} by Perrin-Riou and Kobayashi, respectively, we now define plus/minus subgroups of the local points, as follows.
\begin{definition}[{\cite[Definitions~2.1 and~3.13]{KitOts18}}]\label{def:pm_subgroups} With notations as above we denote, for every $n\geq 1$,
\begin{align*}
\pmformalEC[+](\maxid{n})=\{P\in\formalEC(\maxid{n})\suchthat\trace{n}{m+1}(P)\in\formalEC(\maxid{m})\text{ for all }-1\leq m\leq n-1, m\text{ even}\}
\intertext{and}
\pmformalEC[-](\maxid{n})=\{P\in\formalEC(\maxid{n})\suchthat\trace{n}{m+1}(P)\in\formalEC(\maxid{m})\text{ for all }-1\leq m\leq n-1, m\text{ odd}\}
\end{align*}
Similarly, we set
\begin{align*}	
\pmEC[+](\fullcyc{n})=\{P\in\EC(\fullcyc{n})\suchthat\trace{n}{m+1}(P)\in\EC(\fullcyc{m})\text{ for all }-1\leq m\leq n-1, m\text{ even}\}
\intertext{and}
\pmEC[-](\fullcyc{n})=\{P\in\EC(\fullcyc{n})\suchthat\trace{n}{m+1}(P)\in\EC(\fullcyc{m})\text{ for all }-1\leq m\leq n-1, m\text{ odd}\}
\end{align*}
and we let $\pmEC(\localfield{n})=\pmEC(\fullcyc{n})^\galprime=\HH{0}{\galprime}{\pmEC(\fullcyc{n})}$.
\end{definition}
The next lemma compares the formal signed subgroups of local points with the whole signed subgroups:
\begin{lemma}[{see~\cite[Lemma~3.14]{KitOts18}}] \label{lemma:formal_vs_localPM} Let $\ssprimeid{}=\ssprimeid{i}\in\ssprimeset$ and let $\localfield{}=\localfield[i]{0}$. For all $n\geq 1$ there are exact sequences
\begin{equation}\label{eq:formal_vs_localPM:forPM}
0\longrightarrow\pmformalEC(\maxid{n})\longrightarrow\pmEC(\fullcyc{n})\longrightarrow\pmMod{n}\longrightarrow 0
\end{equation}
where $\pmMod{n}\subseteq\modEC(\finitefield)$ is a finite group, of prime-to-$p$ order bounded independently of $n$. 

More generally, if $\field/\localfield{}$ is any algebraic extension and $\maxid{\field}$ is the maximal ideal of its valuation ring, there is an exact sequence
\begin{equation}
\label{eq:formal_vs_localPM:forbarQ_p}
0\longrightarrow\formalEC(\maxid{\field})\longrightarrow\EC(\field)\longrightarrow D\longrightarrow 0
\end{equation}
where $D$ is a finite group of prime-to-$p$ order, inducing an isomorphism
\[
\tors{\formalEC(\maxid{\overline{\QQ_p}})}{p^\infty}\cong\tors{\EC(\overline{\QQ_p})}{p^\infty}.
\]
\end{lemma}
\begin{proof} Fix $m\geq -1$ and consider the commutative diagram
\[\xymatrix{
0\ar@{->}[0,1]& \formalEC(\maxid{m})\ar@{->}[0,1]\ar@{^(->}[1,0]& \EC(\fullcyc{m})\ar@{->}[0,1]\ar@{^(->}[1,0]& \modEC(\finitefield)\ar@{->}[0,1]\ar@{=}[1,0]& 0\\
0\ar@{->}[0,1]& \formalEC(\maxid{m+1})\ar@{->}[0,1]& \EC(\fullcyc{m+1})\ar@{->}[0,1]& \modEC(\finitefield)\ar@{->}[0,1]& 0
}\]
which induces, by the snake lemma, an isomorphism
\[
\formalEC(\maxid{m+1})/\formalEC(\maxid{m})\overset{\cong}{\longrightarrow} \EC(\fullcyc{m+1})/ \EC(\fullcyc{m}).
\]
Now fix $n\geq m+1$: the above sequence fits into the commutative diagram of exact sequences
\[\xymatrix{
0\ar@{->}[0,1]& \ker\widehat{\trace{n}{m+1}}\ar@{->}[0,1]\ar@{^(->}[1,0]& \ker\overline{\trace{n}{m+1}}\ar@{->}[0,1]\ar@{^(->}[1,0]& \modEC(\finitefield)\ar@{=}[1,0]\\
0\ar@{->}[0,1]& \formalEC(\maxid{n})\ar@{->}[0,1]\ar@{->}[1,0]^{\widehat{\trace{n}{m+1}}}& \EC(\fullcyc{n})\ar@{->}[0,1]\ar@{->}[1,0]^{\overline{\trace{n}{m+1}}}& \modEC(\finitefield)\ar@{->}[0,1]\ar@{->}[1,0]& 0\\
0\ar@{->}[0,1]& \formalEC(\maxid{m+1})/\formalEC(\maxid{m})\ar@{->}[0,1]^{\cong}& \EC(\fullcyc{m+1})/ \EC(\fullcyc{m})\ar@{->}[0,1]& 0\ar@{->}[0,1]& 0
}\]
where $\widehat{\trace{n}{m+1}}$ (\resp $\overline{\trace{n}{m+1}}$) denotes the trace map followed by reduction modulo $\formalEC(\maxid{m+1})$ (\resp modulo $\EC(\fullcyc{m+1})$). In particular, we deduce that $\ker\widehat{\trace{n}{m+1}}$ is a subgroup of $\ker\overline{\trace{n}{m+1}}$ with quotient contained inside $\modEC(\finitefield)$, for every $m\leq n-1$. Taking intersections, we find
\[
\pmformalEC(\maxid{n})=\bigcap_{\substack{(-1)^m=\pm 1\\-1\leq m\leq n-1}}\ker\widehat{\trace{n}{m+1}}\qquad\text{ and }\qquad\pmEC(\fullcyc{n})=\bigcap_{\substack{(-1)^m=\pm 1\\-1\leq m\leq n-1}}\ker\overline{\trace{n}{m+1}}
\]
and therefore an exact sequence
\[
0\longrightarrow\pmformalEC(\maxid{n})\longrightarrow\pmEC(\fullcyc{n})\longrightarrow\pmMod{n}\longrightarrow 0
\]
for some $\pmMod{n}\subseteq \modEC(\finitefield)$. Since $\EC$ has supersingular reduction at $\ssprimeid{}$, the order of $\modEC(\finitefield)$ is prime-to-$p$.

The final isomorphism is simply a translation of the fact that $\modEC(\overline{\finitefield})$ has no $p$-torsion. Using the exact sequence
\[
0\longrightarrow \formalEC(\maxid{\field})\longrightarrow \EC(\field)\longrightarrow \modEC(\rint{\field}/\maxid{\field})\longrightarrow 0,
\]
one obtains 
\[
\tors{\formalEC(\maxid{\field})}{p^\infty}\cong\tors{\EC(\field)}{p^\infty}
\]
and taking direct limit over all $\localfield{}\subseteq\field\subseteq \overline{\QQ_p}$, we deduce the isomorphism in the statement.
\end{proof}
Let $\field\in\{\basefield[n],\fullcyc{n},\genlocal[n]{\genprimeid}\}$ and $G\in \{\GS[n],\Glocal[\fullcyc]{n},\Glocal[{\genlocal[n]}]{\genprimeid}\}$. Recall that for each integer $t\geq 0$, there exists the following functorial exact sequence for $\tors{\EC}{p^t}/\field$
\begin{equation*}
0\longrightarrow \EC(\field)/p^t\EC(\field)\overset{{\fullKum{\field}{p^t}}}{\longrightarrow} \HH{1}{G}{\tors{\EC}{p^t}}\longrightarrow \tors{\HH{1}{G}{\EC}}{p^t}\longrightarrow 0
\end{equation*}
where $\fullKum{\field}{p^t}$ is the Kummer map.
\begin{lemma}[see {\cite[Lemma 8.17]{Kob03}}]\label{lemma:pm_injection} For every $n\geq 1$ and every $t\geq 1$, there is an injection
\[
\pmEC(\localfield{n})/p^t\pmEC(\localfield{n})\longhookrightarrow \EC(\localfield{n})/p^t\EC(\localfield{n})
\]
which induces injections
\[
\pmKum{\localfield{n}}{p^t}\colon \pmEC(\localfield{n})/p^t\pmEC(\localfield{n})\longhookrightarrow \HH{1}{\localfield{n}}{\tors{\EC}{p^t}}.
\]
Similarly, there are injections
\[
\pmKum{\localfield{n}}{p^t}\colon \pmformalEC(\maxid{n}^{\galprime})/p^t\pmformalEC(\maxid{n}^{\galprime})\longhookrightarrow \HH{1}{\localfield{n}}{\tors{\formalEC}{p^t}}.
\]
\end{lemma}
\begin{proof}
Let us first show that
\begin{equation}\label{eq:inj_fullcyclo}
\pmEC(\fullcyc{n})/p^t\pmEC(\fullcyc{n})\longhookrightarrow \EC(\fullcyc{n})/p^t\EC(\fullcyc{n})
\end{equation}
is injective. An element in
\[
\ker\Bigl(\pmEC(\fullcyc{n})/p^t\pmEC(\fullcyc{n})\longrightarrow \EC(\fullcyc{n})/p^t\EC(\fullcyc{n})\Bigr)
\]
is represented by a point $P\in\pmEC(\fullcyc{n})$ such that $P=p^tQ$ for some $Q\in \EC(\fullcyc{n})$. Choose now $m\leq n-1$ such that $(-1)^m=\pm 1$. Taking the trace of $P$ down to $\fullcyc{m+1}$ we obtain that $\trace{n}{m+1}(P)\in\EC(\fullcyc{m})$, by definition of $\pmEC$.
 On the other hand, $\trace{n}{m+1}(P)=p^t\trace{n}{m+1}(Q)$, hence for all $\sigma\in\Gal{\fullcyc{m+1}}{\fullcyc{m}}$ we have $p^t(^\sigma \trace{n}{m+1}(Q)-\trace{n}{m+1}(Q))=0$. Thus $\trace{n}{m+1}(Q)\in\EC(\fullcyc{m})$ thanks to Proposition~\ref{prop:no_p_torsion}, which implies $Q\in\pmEC(\fullcyc{n})$ and the arrow in~\eqref{eq:inj_fullcyclo} is injective. For each $\ast\in\{\emptyset,+,-\}$, taking $\galprime$-cohomology of the tautological exact sequence defining $\pmEC[\ast](\fullcyc{n})/p^t\pmEC[\ast](\fullcyc{n})$ gives
\[
0\longrightarrow \pmEC[\ast](\localfield{n})\overset{\cdot p^t}{\longrightarrow }\pmEC[\ast](\localfield{n})\longrightarrow \HH*{0}{\galprime}{\pmEC[\ast](\fullcyc{n})/p^t\pmEC[\ast](\fullcyc{n})}\longrightarrow\tors{\HH{1}{\galprime}{\pmEC[\ast](\localfield{n})}}{p^t}\longrightarrow 0.
\]
The last module is trivial, because $\galprime$ has order prime-to-$p$, so \begin{equation}\label{eq:galois_invariants_quotient}
\HH*{0}{\galprime}{\pmEC[\ast](\fullcyc{n})/p^t\pmEC[\ast](\fullcyc{n})}=\pmEC[\ast](\localfield{n})/p^t\pmEC[\ast](\localfield{n}).
\end{equation}
Taking $\galprime$-invariants of the injections in~\eqref{eq:inj_fullcyclo} establishes the first part of the lemma. The second is analogous, upon replacing $\pmEC$ with $\pmformalEC$.
\end{proof}
In light of the above Lemma, we can define, for all $n\geq 0$ (and all $\ssprimeid{i}$ if we need to keep track of the local Galois groups), the \emph{signed} Kummer sequence as the exact sequence
\begin{equation}\label{eq:pm_Kummer_sequence}
0\longrightarrow \pmEC(\localfield{n})/p^t\pmEC(\localfield{n})\overset{{\pmKum{\localfield{n}}{p^t}}}{\longrightarrow} \HH{1}{\Glocal{n}}{\tors{\EC}{p^t}}\longrightarrow \HH{1}{\Glocal{n}}{\tors{\EC}{p^t}}/\image{{\pmKum{\localfield{n}}{p^t}}}\longrightarrow 0
\end{equation}
and refer to ${\pmKum{\field}{p^t}}$ as the \emph{signed} Kummer map. Analogous signed Kummer exact sequence can be defined for the formal group $\formalEC$, as follows:
\begin{equation}\label{eq:pm_Kummer_formal}
0\longrightarrow \pmformalEC(\localfield{n})/p^t\pmformalEC(\maxid{n}^{\galprime})\overset{{\pmKum{\localfield{n}}{p^t}}}{\longrightarrow} \HH{1}{\Glocal{n}}{\tors{\formalEC}{p^t}}\longrightarrow \HH{1}{\Glocal{n}}{\tors{\formalEC}{p^t}}/\pmKum{\localfield{n}}{p^t}\bigl(\pmformalEC(\maxid{n}^{\galprime})\bigr)\longrightarrow 0
\end{equation}
\begin{remark}
It is perhaps interesting to stress that the signed Kummer map defined in~\eqref{eq:pm_Kummer_sequence} \emph{does not arise} as a connecting homomorphism in Galois cohomology. Indeed, $\pmEC$ is only defined at the level of points for extensions in the cyclotomic tower and it is not a sub-representation of $\EC$, since in the supersingular case the local Galois representation $\tors{\EC}{p^\infty}$ is irreducible.
\end{remark}
\subsection{The multi-signed Selmer groups}
We use the notation introduced in~\ref{not:S_and_primes}. For generalities regarding the classical Selmer group for $\tors{\EC}{p^t}/\basefield[n]$ (for $1\leq t<\infty$) and for $\tors{\EC}{p^\infty}/\basefield[n]$ we refer to \cite[Chapters~1 and~2]{CoaSuj10}. They are defined as
\begin{equation*}
\begin{split}
 \fullSel{\tors{\EC}{p^t}}{\basefield[n]}&=\ker\Bigl(\HH{1}{\GS[n]}{\tors{\EC}{p^t}}\longrightarrow \bigoplus_{\genprimeid\in\primeset[{\basefield[n]}]}\tors{\HH{1}{\genlocal[n]{\genprimeid}}{\EC}}{p^t}\Bigr)\\
&=\ker\Bigl(\HH{1}{\GS[n]}{\tors{\EC}{p^t}}\longrightarrow \bigoplus_{\genprimeid\in\primeset[{\basefield[n]}]\setminus\ssprimeset[{\basefield[n]}]}\tors{\HH{1}{\genlocal[n]{\genprimeid}}{\EC}}{p^t}\oplus\bigoplus_{i=1}^{\ordssprimeset}\tors{\HH{1}{\localfield[i]{n}}{\EC}}{p^t}\Bigr)\\
&=\ker\Bigl(\HH{1}{\GS[n]}{\tors{\EC}{p^t}}\longrightarrow \bigoplus_{\genprimeid\in\primeset[{\basefield[n]}]\setminus\ssprimeset[{\basefield[n]}]}\HH{1}{\genlocal[n]{\genprimeid}}{\tors{\EC}{p^t}}/\image\fullKum{\localfield[i]{n}}{p^t}\oplus\bigoplus_{i=1}^{\ordssprimeset}\HH{1}{\localfield[i]{n}}{\tors{\EC}{p^t}}/\image\fullKum{\localfield[i]{n}}{p^t} \Bigr)
\end{split}\end{equation*}
where we split the sum into two parts, one over primes in $\ssprimeset$ and the over primes which are not supersingular. This is done mainly for future comparison with signed Selmer groups.
Passing to the limit over $t$, one defines
\begin{equation*}
 \fullSel{\tors{\EC}{p^\infty}}{\basefield[n]}=\varinjlim_t\fullSel{\tors{\EC}{p^t}}{\basefield[n]}.
\end{equation*}
In the supersingular reduction case, the Iwasawa theory of the signed Selmer groups as initially defined by Perrin-Riou and Kobayashi respectively in~\cite{Per90} and~\cite{Kob03} is of particular interest. The multi-signed residual Selmer groups are defined below and we postpone a larger discussion, from the Iwasawa-theoretic point of view, to Section~\ref{section:iwasawa}. Our main reference is the work~\cite{KitOts18} by Kitajima--Otsuki. Fix a vector of signs $\multisign\in\{+,-\}^{\ordssprimeset}$ throughout. For every local place $\ssprimeid{i}$ with $1\leq i\leq\ordssprimeset$, to ease notation write $\pmKum{\localfield[i]{n}}{p^t}$ to denote either $\pmKum[+]{\localfield[i]{n}}{p^t}$ or $\pmKum[-]{\localfield[i]{n}}{p^t}$, according to the $i$th component $\multisign_i$ of the vector of signs.
\begin{definition}\label{def:pm-Sel}
For every intermediate number field $\basefield\subseteq \basefield[n]\subsetneq \basecyc$ define the \emph{fine multi-signed residual Selmer group} as
\begin{align*}
\pmSel{\tors{\EC}{p}}{\basefield[n]}=\ker\Bigl(\HH{1}{\GS[n]}{\tors{\EC}{p}}\rightarrow \bigoplus_{\badprimeid{\empty}\in\badprimeset[{\basefield[n]}]}\HH{1}{\genlocal[n]{\badprimeid{\empty}}}{\tors{\EC}{p}}\oplus\bigoplus_{\ordprimeid{\empty}\in\ordprimeset[{\basefield[n]}]}\HH{1}{\genlocal[n]{\ordprimeid{\empty}}}{\tors{\modEC}{p}}\oplus\bigoplus_{i=1}^{\ordssprimeset}\HH{1}{\localfield[i]{n}}{\tors{\EC}{p}}/\image \pmKum{\localfield[i]{n}}{p}\Bigr),
\end{align*}
where, at an ordinary prime $\ordprimeid{\empty}$, $\tors{\modEC}{p}$ is seen as a $\Glocal[{\genlocal[n]}]{\ordprimeid{\empty}}$-module through the surjection $\Glocal[{\genlocal[n]}]{\ordprimeid{\empty}}\twoheadrightarrow \Glocal[{\genlocal[n]}]{\ordprimeid{\empty}}^\mathrm{ur}=\Glocal[{\finitefield[\ordprimeid{\empty}]}]{\empty}$. Similarly, the \emph{usual multi-signed Selmer group} is defined as
\begin{align*}
\KobSel{\tors{\EC}{p^\infty}}{\basefield[n]}=\ker\Bigl(\HH{1}{\GS[n]}{\tors{\EC}{p^\infty}}\rightarrow \bigoplus_{\genprimeid\in\primeset[{\basefield[n]}]\setminus\ssprimeset[{\basefield[n]}]}\tors{\HH{1}{\genlocal[n]{\genprimeid}}{\EC}}{p^\infty}\oplus\bigoplus_{i=1}^{\ordssprimeset}\HH{1}{\localfield[i]{n}}{\tors{\EC}{p^\infty}}/\image \pmKum{\localfield[i]{n}}{p^\infty}\Bigr).
\end{align*}
\end{definition}
As the notation suggests, these multi-signed residual Selmer groups only depend upon the isomorphism class of $\tors{\EC}{p}$ rather than on the curve $\EC$ itself, at least when assuming~\ref{hyp_1}. This is the content of Corollary~\ref{cor:Sel_indep}, which relies on Proposition~\ref{prop:universal_sequence} and Proposition~\ref{prop:raynaud} below.

We start with the following technical lemma:
\begin{lemma}\label{lemma:H[p]=H(p)_finite} Fix $n\geq 0$ and let $G\in\{\GS[n],\Glocal{n},\Glocal[{\genlocal[n]}]{\badprimeid{}}\}$ for some $\badprimeid{}\in \badprimeset$. Denote by $\fakeKum[G]{n}=\fakeKum{n}$ the natural surjective arrow 
\[
\fakeKum{n}\colon\HH{1}{G}{\tors{\EC}{p}}\longrightarrow \tors{\HH{1}{G}{\tors{\EC}{p^\infty}}}{p}.
\]
For $\ordprimeid{}\in\ordprimeset$ and $G=\Glocal[{\genlocal[n]}]{\ordprimeid{\empty}}$  let
\[
\modfakeKum[G]{n}=\modfakeKum{n}\colon\HH{1}{G}{\tors{\modEC}{p}}\longrightarrow \tors{\HH{1}{G}{\tors{\modEC}{p^\infty}}}{p}
\]
be the analogous surjection for the Galois representation $\tors{\modEC}{p^\infty}$. Then the following assertions hold.
\begin{enumerate}[label=\roman*)]
\item\label{point:lemma:H[p]=H(p)_finite_ss:global} If $G\in\{\GS[n],\Glocal{n}\}$, then $\fakeKum{n}$ is an isomorphism.
\item\label{point:lemma:H[p]=H(p)_finite_bad} If $G=\Glocal[{\genlocal[n]}]{\badprimeid{}}$ for some $\badprimeid{}\in \badprimeset$, then $\ker{\fakeKum{n}}$ is an $\finitefield$-vector space of dimension $\dim_{\finitefield}(\ker{\fakeKum{n}})=\dim_{\finitefield}\tors{\EC(\genlocal[n]{\badprimeid{}})}{p}\leq 2$.
\item\label{point:lemma:H[p]=H(p)_finite_ord} If $G=\Glocal[{\genlocal[n]}]{\ordprimeid{\empty}}$ for some $\ordprimeid{\empty}\in \ordprimeset$, then $\ker{\modfakeKum{n}}$ is an $\finitefield$-vector space of dimension $\dim_{\finitefield}(\ker{\modfakeKum{n}})=\dim_{\finitefield}\tors{\modEC(\genlocal[n]{\ordprimeid{\empty}})}{p}\leq 1$.
\end{enumerate}
\end{lemma}
\begin{proof}
Taking $G$-cohomology of the exact sequence
\begin{equation}\label{eq:false_kummer}
0\longrightarrow \tors{\EC}{p}\longrightarrow \tors{\EC}{p^\infty}\longrightarrow \tors{\EC}{p^\infty}\longrightarrow 0
\end{equation}
gives an exact sequence
\begin{equation}\label{eq:cohom_pinfty=p}
0\longrightarrow \ker{\fakeKum{n}}=\HH{0}{G}{\tors{\EC}{p^\infty}}/p\HH{0}{G}{\tors{\EC}{p^\infty}}\longrightarrow \HH{1}{G}{\tors{\EC}{p}}\overset{\fakeKum{n}}{\longrightarrow} \tors{\HH{1}{G}{\tors{\EC}{p^\infty}}}{p}\longrightarrow 0.
\end{equation}

When $G\in\{\GS[n],\Glocal{n}\}$, the first term in~\eqref{eq:cohom_pinfty=p} is trivial thanks to Proposition~\ref{prop:no_p_torsion} and assertion~\ref{point:lemma:H[p]=H(p)_finite_bad} follows.

When $G=\Glocal[{\genlocal[n]}]{\badprimeid{}}$ for some $\badprimeid{}\in \badprimeset$, the first term in~\eqref{eq:cohom_pinfty=p} has $\finitefield$-dimension equal to $\dim_{\finitefield}\tors{\EC(\genlocal[n]{\badprimeid{}})}{p}$, since $\tors{\EC(\genlocal[n]{\badprimeid{}})}{p^\infty}$ is finite. Moreover, the group $\HH{0}{G}{\tors{\EC}{p}}$ is a subgroup of $\tors{\EC(\overline{\genlocal{\badprimeid{}}})}{p}\cong(\finitefield)^2$. This shows that this dimension is bounded by $2$, whence assertion~\ref{point:lemma:H[p]=H(p)_finite_bad}.

Finally, when $G=\Glocal[{\genlocal[n]}]{\ordprimeid{\empty}}$ for some $\ordprimeid{\empty}\in \ordprimeset$, replace~\eqref{eq:false_kummer} by
\begin{equation*}
0\longrightarrow \tors{\modEC}{p}\longrightarrow \tors{\modEC}{p^\infty}\longrightarrow \tors{\modEC}{p^\infty}\longrightarrow 0
\end{equation*}
to obtain an exact sequence
\begin{equation}\label{eq:mod_cohom_pinfty=p}
0\longrightarrow \ker{\modfakeKum{n}}=\HH{0}{G}{\tors{\modEC}{p^\infty}}/p\HH{0}{G}{\tors{\modEC}{p^\infty}}\longrightarrow \HH{1}{G}{\tors{\modEC}{p}}\overset{\modfakeKum{n}}{\longrightarrow} \tors{\HH{1}{G}{\tors{\modEC}{p^\infty}}}{p}\longrightarrow 0.
\end{equation}
The first term in~\eqref{eq:mod_cohom_pinfty=p} is a $\finitefield$-vector space of dimension bounded by $\dim_{\finitefield}\tors{\modEC(\overline{\finitefield[\ordprimeid{\empty}]})}{p^\infty}/p\tors{\modEC(\overline{\finitefield[\ordprimeid{\empty}]})}{p^\infty}\cong\finitefield$. This finishes the proof.
\end{proof}

Let us now prove that the local conditions in the definition of the multi-signed residual Selmer group depend only on the residual representation also for primes above $p$, beginning with supersingular primes. Under our standing assumption~\ref{hyp_1}, $\EC$ is supersingular at all primes $\ssprimeid{1},\hdots,\ssprimeid{\ordssprimeset}$ and the exact sequence~\eqref{eq:formal_vs_localPM:forbarQ_p} induces an isomorphism $\HH{1}{\localfield[i]{n}}{\tors{\EC}{p^\infty}}\cong\HH{1}{\localfield[i]{n}}{\tors{\formalEC}{p^\infty}}$, for all $n\geq 0$. On the other hand, the exact sequence~\eqref{eq:formal_vs_localPM:forPM} shows that the images of the signed Kummer maps $\pmKum{\localfield[i]{n}}{p^\infty}\bigl(\pmEC\bigl(\localfield{i,n})\bigr)$ and $\pmKum{\localfield[i]{n}}{p^\infty}\bigl(\pmformalEC\bigl(\maxid{n}^{\galprime})\bigr)$ are isomorphic. It is straightforward to check that these isomorphisms are compatible, and in turn induce isomorphisms
\begin{equation}\label{eq:iso_H1_EC_formal}
\HH{1}{\localfield[i]{n}}{\tors{\EC}{p^\infty}}/\image(\pmKum{\localfield[i]{n}}{p^\infty})\cong \HH{1}{\localfield[i]{n}}{\tors{\formalEC}{p^\infty}}/\pmKum{\localfield[i]{n}}{p^\infty}\bigl(\pmformalEC(\maxid{n}^{\galprime})\bigr).
\end{equation}
As discussed in \cite[Corollary~8.5]{Kob03} and \cite[\S~3.1]{KitOts18}, there is a $\rint{\localfield{0}}=\ZZ_p$-isomorphism 
\begin{equation}\label{eq:iso_Fss_formal}
\log_{\Fss}\circ \exp_{\formalEC}\colon\formalEC\overset{\cong}{\longrightarrow}\Fss
\end{equation}
where $\Fss$ is the supersingular formal group whose logarithm $\log_{\Fss}$ is of Honda type $t^2+p$. In particular, the isomorphism class of the formal group $\formalEC$ is independent of the curve $\EC$, whenever the curve satisfies~\ref{hyp_1}. Moreover, for every $n$ there are two subgroups $\pmFss(\maxid{n})\subseteq \Fss(\maxid{n})$ defined by the same norm relations defining $\pmformalEC$ (see Definition~\ref{def:pm_subgroups}), but for points on the formal group $\Fss$ rather than $\formalEC$. Equivalently, they are defined as
\[
\pmFss(\maxid{n})= (\log_{\Fss}\circ\exp_{\formalEC})\bigl(\pmformalEC(\maxid{n})\bigr).
\]
Therefore, as subgroups of $\Fss(\maxid{n})$, they are independent of $\EC$ for all $n\geq 0$. Moreover, there is an evident definition of the analogues of the signed Kummer sequence~\eqref{eq:pm_Kummer_formal} for $\Fss$ and $\pmFss$ instead of $\formalEC$. Combining~\eqref{eq:iso_H1_EC_formal} with~\eqref{eq:iso_Fss_formal} gives
\[
\HH{1}{\localfield[i]{n}}{\tors{\EC}{p^\infty}}/\image(\pmKum{\localfield[i]{n}}{p^\infty})\cong\HH{1}{\localfield[i]{n}}{\tors{(\Fss)}{p^\infty}}/\pmKum{\localfield[i]{n}}{p^\infty}\bigl(\pmFss(\maxid{n}^{\galprime})\bigr),
\]
where the right-hand side does not depend on $\EC$. We summarise the above discussion in the following
\begin{proposition}\label{prop:universal_sequence} Let $\EC/\basefield$ be an elliptic curve satisfying hypothesis~\ref{hyp_1}. For all $1\leq i\leq \ordssprimeset$ and all $n\geq 0$, there are functorial isomorphisms
\begin{equation*}
\HH{1}{\localfield[i]{n}}{\tors{\EC}{p^\infty}}/\image(\pmKum{\localfield[i]{n}}{p^\infty})\cong \HH{1}{\localfield[i]{n}}{\tors{(\Fss)}{p^\infty}}/\pmKum{\localfield[i]{n}}{p^\infty}\bigl(\pmFss(\maxid{n}^{\galprime})\bigr)
\end{equation*}
In particular, the modules $\HH{1}{\localfield[i]{n}}{\tors{\EC}{p^\infty}}/\image(\pmKum{\localfield[i]{n}}{p^\infty})$ are independent of $\EC$, since the right-hand sides are.
\end{proposition}
\begin{proof} The fact that the first isomorphism is functorial follows from Honda theory, which shows that the isomorphism between $\Fss$ and $\formalEC$ is given by $\log_{\formalEC}\circ \exp_{\Fss}$ (see~\cite[Theorem~8.3~(ii)]{Kob03}).
\end{proof}

What remains to be proven is the analogue of the above result when replacing $\tors{\EC}{p^\infty}$ by $\tors{\EC}{p}$, which is the module we are ultimately interested in. This is done in Proposition~\ref{prop:H[p]=H(p)_infinite}-\ref{point:prop:H[p]=H(p)_infinite_ss}, as we move up the cyclotomic tower. Concerning ordinary primes, we have the following result.
\begin{proposition}\label{prop:raynaud} Let $\EC_1,\EC_2$ be two elliptic curves defined over $\basefield$ satisfying hypothesis~\ref{hyp_1}. Let $\ssprimeset[j]$ \textup{(}\resp $\ordprimeset[j]$\textup{)} denote the set of primes where $\EC_j$ has supersingular \textup{(}\resp ordinary\textup{)} reduction, for $j=1,2$. Then
\begin{enumerate}[label=\roman*)]
\item\label{item:prop:raynaud_Sss=Sord} $\ssprimeset[1]=\ssprimeset[2]$ and $\ordprimeset[1]=\ordprimeset[2]$. Denote these sets simply by $\ssprimeset$ and $\ordprimeset$, respectively.
\item\label{item:prop:raynaud_modE=modE} Every isomorphism $\tors{(\EC_1)}{p}\cong\tors{(\EC_2)}{p}$ induces an isomorphism $\tors{(\modEC_1)}{p}\cong\tors{(\modEC_2)}{p}$ and, in particular, an isomorphism
\[
\HH{1}{\genlocal[n]{\ordprimeid{\empty}}}{\tors{(\modEC_1)}{p}}\cong \HH{1}{\genlocal[n]{\ordprimeid{\empty}}}{\tors{(\modEC_2)}{p}}\qquad\text{ for all }n\geq 0.
\]
\end{enumerate}
\end{proposition}
\begin{proof}
Starting with~\ref{item:prop:raynaud_Sss=Sord}, observe that an equality $\ssprimeset[1]=\ssprimeset[2]$ will imply $\ordprimeset[1]=\ordprimeset[2]$ because $\ordprimeset[j]=\pprimeset\setminus\ssprimeset[j]$ (by~\ref{hyp_1}, both curves have good reduction at all primes in $\pprimeset$). To show the claimed equality, pick a prime $\ssprimeid{}\in\ssprimeset[1]$. By~\ref{hyp_1}, $\genlocal{\ssprimeid{\empty}}\cong \QQ_p$ is absolutely unramified. Denote by $\neron_j$ the Néron model of $\EC_j$. (see~\cite[IV,~Corollary~6.3]{Sil94} for the existence of this model). Note that the operations of passing to the generic (\resp special) fibre and of computing the kernel of multiplication by $p$ are fibre products. Thus these two operations commute and, in particular, the generic (\resp special) fibre of the finite, flat $\rint{\genlocal{\ssprimeid{\empty}}}$-group scheme $\tors{(\neron_j)}{p}$ is isomorphic to ${\tors{(\EC_j)}{p}}$ (\resp to $\tors{(\modEC_j)}{p}$), for $j=1,2$. Applying~\cite[Corollaire~3.3.6]{Ray74}, we see that the hypothesis $\tors{(\EC_1)}{p}\cong\tors{(\EC_2)}{p}$ (as Galois modules or, what amounts to the same, as finite, flat $\genlocal{\ssprimeid{}}$-group schemes) grants the existence of an isomorphism
\begin{equation}\label{eq:iso_raynaud}
\tors{(\neron_1)}{p}\cong\tors{(\neron_2)}{p}
\end{equation}
of finite, flat $\rint{\genlocal{\ssprimeid{}}}$-group schemes. By taking closed fibres, this yields an isomorphism
\[
{\tors{(\modEC_1)}{p}}\cong{\tors{(\neron_1)}{p}}_{/\finitefield[\ssprimeid{}]}\cong{\tors{(\neron_2)}{p}}_{/\finitefield[\ssprimeid{}]}\cong{\tors{(\modEC_2)}{p}}
\]
as finite flat $\finitefield[\ssprimeid{}]$-group schemes. This shows that the elliptic cuve $\modEC_2$ has supersingular reduction at $\ssprimeid{}$ and $\ssprimeset[1]\subseteq\ssprimeset[2]$. By reversing the role of $\EC_1$ and $\EC_2$, this yields $\ssprimeset[1]=\ssprimeset[2]$, as claimed.

Passing to~\ref{item:prop:raynaud_modE=modE}, let $\ordprimeid{}$ be a prime where one, and hence both curves, have ordinary reduction. The assumption $e(\ordprimeid{})<p-1$ allows us to again apply~\cite[Corollaire~3.3.6]{Ray74} and the isomorphism~\eqref{eq:iso_raynaud} holds again, where now $\neron_j$ denotes the Néron model of $\EC_j/\genlocal{\ordprimeid{}}$. Taking closed fibres, we obtain
\[
{\tors{(\modEC_1)}{p}}\cong{\tors{(\neron_1)}{p}}_{/\finitefield[\ordprimeid{}]}\cong{\tors{(\neron_2)}{p}}_{/\finitefield[\ordprimeid{}]}\cong{\tors{(\modEC_2)}{p}}
\]
finishing the proof of the proposition.
\end{proof}

\section{Iwasawa theory for the signed Selmer groups}\label{section:iwasawa}
\subsection{Cyclotomic multi-signed residual Selmer groups}
In this section we focus on the Iwasawa theory for the multi-signed residual Selmer group introduced in Definition~\ref{def:pm-Sel}. Retaining the notation introduced in~\ref{not:S_and_primes}, set $\GS[\cyc]=\Gal{\baseS}{\basecyc}$. Denote by $\Gamma$ the Galois group $\Gal{\basecyc}{\basefield}\cong\Gal{\ssloccyc}{\QQ_p}$, let $\IwLambda=\ZZ_p[\![\Gamma]\!]$ be its Iwasawa algebra, and set $\resLambda=\finitefield[p][\![\Gamma]\!]$. For any module $M$ over an Iwasawa algebra, its Pontryagin dual $\Hom_{\ZZ_p}(M,\QQ_p/\ZZ_p)$ is denoted by $\pontr{M}$. When $M$ is discrete, we say that it is cofinitely generated (\resp cofree, cotorsion, of corank equal to $m\in\NN$) to mean that $\pontr{M}$ is finitely generated (\resp free, torsion, of rank equal to $m$) over the Iwasawa algebra. Observe that, given any co-finitely generated $\IwLambda$-module $M$, there is an equality $\pontr{M}/p\pontr{M}=\pontr{(\tors{M}{p})}$, inducing the inequality
\[
\corank M\leq \corank[\resLambda]\tors{M}{p}
\]
which is an equality if and only if the $\iwmu[]$ invariant of $\pontr{M}$ vanishes.

Thanks to Lemma~\ref{lemma:pm_injection}, there is an inclusion of subgroups in $\HH{1}{\localfield{n}}{\tors{\EC}{p}}$
\[
\image(\pmKum{\localfield{n}}{p})\longhookrightarrow \image(\fullKum{\localfield{n}}{p}),
\]
which will play a role in defining the Selmer groups. We display the subscript $1\leq i\leq \ordssprimeset$ to keep track of the local Galois cohomology groups, writing $\image(\pmKum{\localfield[i]{n}}{p})\hookrightarrow \image(\fullKum{\localfield[i]{n}}{p})\subseteq\HH{1}{\localfield[i]{n}}{\tors{\EC}{p}}$.

By taking the direct limit of the exact sequence~\eqref{eq:pm_Kummer_sequence} over the subextensions inside $\localfield{\cyc}/\localfield{}$ gives the exact sequences, for all $1\leq t\leq \infty$,
\begin{equation}\label{eq:pm_Kummer_sequence_cyclo}
0\longrightarrow \pmEC(\localfield{\cyc})/p^t\pmEC(\localfield{\cyc})\overset{{\pmKum{\localfield{\cyc}}{p^t}}}{\longrightarrow} \HH{1}{\localfield{\cyc}}{\tors{\EC}{p^t}}\longrightarrow \HH{1}{\localfield{\cyc}}{\EC}/\image(\pmKum{\localfield[i]{\cyc}}{p^t})\longrightarrow 0,
\end{equation}
The following proposition is the main technical tool needed to compare local and global cohomology groups of the residual representation $\tors{\EC}{p}$ along the cyclotomic tower, with those of the representation $\tors{\EC}{p^\infty}$. We refer to Lemma~\ref{lemma:H[p]=H(p)_finite} for the definition of the arrows $\fakeKum{\empty}$ in the statement below.
\begin{proposition}\label{prop:H[p]=H(p)_infinite} Let $G\in\{\GS[\cyc],\Glocal{\cyc},\Glocal[{\genlocal[\cyc]}]{\GenPrimeId}\}$ where $\GenPrimeId\mid\genprimeid\in \badprimeset\cup\ordprimeset$. Write $\pmKum{\genlocal[\cyc]{\GenPrimeId}}{p^\infty}$ to denote $\fullKum{\genlocal[\cyc]{\GenPrimeId}}{p^\infty}$ when $\GenPrimeId\mid\genprimeid\in \badprimeset\cup\ordprimeset$ \textup{(}in particular, these maps are independent of the sign $\pm$\textup{)}.
\begin{enumerate}[label=\alph*\textup{)}]
\item \label{point:prop:H[p]=H(p)_infinite_global}If $G\in\{\GS[\cyc],\Glocal{\cyc}\}$, the map $\fakeKum[G]{\cyc}$ is an isomorphism $\HH{1}{G}{\tors{\EC}{p}}\xrightarrow{\cong} \tors{\HH{1}{G}{\tors{\EC}{p^\infty}}}{p}$.
\item \label{point:prop:H[p]=H(p)_infinite_bad} If $G=\Glocal[{\genlocal[\cyc]}]{\GenPrimeId}$ for some $\GenPrimeId\mid\badprimeid{\empty}\in \badprimeset$, the kernel of $\fakeKum[\GenPrimeId]{\cyc}\colon\HH{1}{G}{\tors{\EC}{p}}\twoheadrightarrow\tors{\HH{1}{G}{\tors{\EC}{p^\infty}}}{p}$ is finite, of dimension $\dim_{\finitefield}(\ker{\fakeKum[\GenPrimeId]{\cyc}})=\dim_{\finitefield}\tors{\EC(\genlocal[\cyc]{\GenPrimeId})}{p}\leq 2$, and
\begin{equation*}
\corank[\resLambda]\HH{1}{G}{\tors{\EC}{p}}=\corank[\resLambda]\tors{\HH{1}{G}{\tors{\EC}{p^\infty}}}{p}.
\end{equation*}
\item \label{point:prop:H[p]=H(p)_infinite_ord}If $G=\Glocal[{\genlocal[\cyc]}]{\GenPrimeId}$ for some $\GenPrimeId\mid\ordprimeid{\empty}\in \ordprimeset$, $\modfakeKum[\GenPrimeId]{\cyc}$ extends to a surjective map
\[
\modfakeKum[\GenPrimeId]{\cyc}\colon\HH{1}{G}{\tors{\modEC}{p}}\longtwoheadrightarrow\tors{\HH{1}{G}{\EC}}{p}=\tors{\left(\HH{1}{G}{\tors{\EC}{p^\infty}}/\image(\pmKum{\genlocal[\GenPrimeId]{\cyc}}{p^\infty})\right)}{p}
\]
whose kernel is finite, of dimension $\dim_{\finitefield}(\ker{\modfakeKum[\GenPrimeId]{\cyc}})=\dim_{\finitefield}\tors{\modEC(\finitefield[\GenPrimeId])}{p}\leq 1$, and
\begin{equation*}
\corank[\resLambda]\HH{1}{G}{\tors{\modEC}{p}}=\corank[\resLambda]\tors{\left(\HH{1}{G}{\tors{\EC}{p^\infty}}/\image(\pmKum{\genlocal[\GenPrimeId]{\cyc}}{p^\infty})\right)}{p}.
\end{equation*}
\item \label{point:prop:H[p]=H(p)_infinite_ss} If $G=\Glocal{\cyc}$, the morphism $\fakeKum[\ssprimeid]{\cyc}$ induces an isomorphism
\begin{equation*}
\pmfakeKum[\ssprimeid]{\cyc}\colon\HH{1}{G}{\tors{\EC}{p}}/\image(\pmKum{\ssloccyc}{p})\overset{\cong}{\longrightarrow}\tors{\left(\HH{1}{G}{\tors{\EC}{p^\infty}}/\image(\pmKum{\ssloccyc}{p^\infty})\right)}{p}
\end{equation*}
giving
\begin{equation*}
\corank[\resLambda]\Bigl(\HH{1}{G}{\tors{\EC}{p}}/\image(\pmKum{\ssloccyc}{p})\Bigr)=\corank[\resLambda]\tors{\left(\HH{1}{G}{\tors{\EC}{p^\infty}}/\image(\pmKum{\ssloccyc}{p^\infty})\right)}{p}.
\end{equation*}
\end{enumerate}
\end{proposition}
\begin{proof}
The isomorphism in~\ref{point:prop:H[p]=H(p)_infinite_global} follows immediately from passing to the direct limit of the isomorphisms at finite levels, proven in~Lemma~\ref{lemma:H[p]=H(p)_finite}-\ref{point:lemma:H[p]=H(p)_finite_ss:global}. Similarly, the description of the kernels in~\ref{point:prop:H[p]=H(p)_infinite_bad} follows from passing to the direct limit in~Lemma~\ref{lemma:H[p]=H(p)_finite}-\ref{point:lemma:H[p]=H(p)_finite_bad}.

The proof of~\ref{point:prop:H[p]=H(p)_infinite_ord} relies on the theory of \emph{deeply ramified extensions} as defined by Coates and Greenberg (see~\cite{CoaGre96}, in particular Theorem~2.13 \ibid,~noting that the cyclotomic $\ZZ_p$-extension is deeply ramified). 
Consider the exact sequence
\[
0\longrightarrow \tors{\ordformalEC}{p^\infty}\longrightarrow\tors{\EC}{p^\infty}\longrightarrow\tors{\modEC}{p^\infty}\longrightarrow 0,
\]
where $\ordformalEC$ is the formal group of $\EC/\rint{\genlocal[\cyc]{\ordprimeid{\empty}}}$. Then the long exact $G$-cohomology sequence gives
\begin{equation}\label{eq:CG_formalEC}
0\longrightarrow\HH{1}{G}{\tors{\EC}{p^\infty}}/\image(\HH{1}{G}{\tors{\ordformalEC}{p^\infty}})\longrightarrow\HH{1}{G}{\tors{\modEC}{p^\infty}}\longrightarrow\HH{2}{G}{\tors{\ordformalEC}{p^\infty}}.
\end{equation}
We claim that $\HH{2}{G}{\tors{\ordformalEC}{p^\infty}}=0$. Indeed, $\HH{2}{G}{\tors{\ordformalEC}{p^\infty}}=\varinjlim \HH{2}{G}{\tors{\ordformalEC}{p^t}}$ and it will be enough to show that $\HH{2}{G}{\tors{\ordformalEC}{p^t}}=0$ for all $t\geq 0$. This follows from the fact that $G$ has $p$-cohomological dimension $1$ (see~\cite[proof~of~Proposition~9,~Chapitre~II,~\S3.3]{Ser94}).

Thus we obtain from~\eqref{eq:CG_formalEC} an isomorphism
\[
\HH{1}{G}{\tors{\EC}{p^\infty}}/\image(\HH{1}{G}{\tors{\ordformalEC}{p^\infty}})\cong\HH{1}{G}{\tors{\modEC}{p^\infty}}.
\]
By \cite[Proposition~4.3 and diagram~(4.8)]{CoaGre96}, we further have
\begin{equation}\label{eq:CoaGre:modp}
\tors{\HH{1}{G}{\EC}}{p^\infty}=\HH{1}{G}{\tors{\EC}{p^\infty}}/\image(\pmKum{\genlocal[\cyc]{\GenPrimeId}}{p^\infty})=\HH{1}{G}{\tors{\EC}{p^\infty}}/\image(\HH{1}{G}{\tors{\ordformalEC}{p^\infty}})\cong\HH{1}{G}{\tors{\modEC}{p^\infty}}.
\end{equation}
Hence $\tors{\HH{1}{G}{\EC}}{p^\infty}\cong\HH{1}{G}{\tors{\modEC}{p^\infty}}$ and, in particular,
\begin{equation*}
\tors{\HH{1}{G}{\EC}}{p}\cong\tors{\HH{1}{G}{\tors{\modEC}{p^\infty}}}{p}.
\end{equation*}
It follows that the surjective arrow $\modfakeKum[\GenPrimeId]{\cyc}\colon\HH{1}{G}{\tors{\modEC}{p}}\twoheadrightarrow\tors{\HH{1}{G}{\tors{\modEC}{p^\infty}}}{p}$ takes values in $\tors{\HH{1}{G}{\EC}}{p}$ and its kernel is finite, of $\finitefield$-dimension less or equal to~$1$, by Lemma~\ref{lemma:H[p]=H(p)_finite}-\ref{point:lemma:H[p]=H(p)_finite_ord}.

The equality of $\resLambda$-coranks in~\ref{point:prop:H[p]=H(p)_infinite_bad} and~\ref{point:prop:H[p]=H(p)_infinite_ord} follows from the fact that a finite module has trivial $\resLambda$-rank.

To prove assertion~\ref{point:prop:H[p]=H(p)_infinite_ss}, note that $\pmEC(\localfield{\cyc})$ is $p$-torsion free. Indeed, $\tors{\pmEC(\localfield{\cyc})}{p}=\tors{\pmformalEC(\maxid{\infty}^{\galprime})}{p}$ by Lemma~\ref{lemma:formal_vs_localPM}. But $\tors{\pmformalEC(\maxid{\cyc})}{p}\subseteq \tors{\formalEC(\maxid{\cyc})}{p}=0$, by Proposition~\ref{prop:no_p_torsion}, hence $\tors{\pmEC(\localfield{\cyc})}{p}=0$. In particular, $\pmEC(\localfield{\cyc})$ is a direct limit of free $\ZZ_p$-modules of finite rank, hence
$\operatorname{Tor}^1_{\ZZ_p}(\pmEC(\localfield{\cyc}),\QQ_p/\ZZ_p)=0$. Consider the exact sequence
\begin{align*}
0\longrightarrow \ZZ/p\longrightarrow \QQ_p/\ZZ_p\longrightarrow \QQ_p/\ZZ_p\longrightarrow 0.
\end{align*}
Tensoring it with $\pmEC(\localfield{\cyc})$ over $\ZZ_p$ yields
\begin{equation}\label{eq:EC/p=EC[p]}
\pmEC(\localfield{\cyc})\otimes\ZZ/p\cong\tors{(\pmEC(\localfield{\cyc})\otimes\QQ_p/\ZZ_p)}{p}.
\end{equation}
Now, the map $\pmfakeKum[\ssprimeid]{\cyc}$, defined as the composition of $\fakeKum[\ssprimeid]{\cyc}$ with reduction modulo $\image\pmKum{\localfield{\cyc}}{p^\infty}$, appears in the following diagram of exact sequences:
\begin{equation}\begin{split}\label{diag:Kumm_modp_ss}
\xymatrix@C=2.5em{
0\ar@{->}[0,1]&\pmEC(\localfield{\cyc})\otimes\ZZ/p\ar@{->}[0,1]^-{\pmKum{\localfield{\cyc}}{p}}\ar@{->}[1,0]_{\rotatebox{270}{$\cong$}}&\HH{1}{\localfield{\cyc}}{\tors{\EC}{p}}\ar@{->}[0,1]\ar@{->}[1,0]_{\rotatebox{270}{$\cong$}}^{\fakeKum{\GS[cyc]}}&\HH{1}{\localfield{\cyc}}{\tors{\EC}{p}}/\image\pmKum{\localfield{\cyc}}{p}\ar@{->}[0,1]\ar@{->}[1,0]^{\pmfakeKum[\ssprimeid]{\cyc}}& 0\\
0\ar@{->}[0,1]&\tors{\Bigl(\pmEC(\localfield{\cyc})\otimes\QQ_p/\ZZ_p\Bigr)}{p}\ar@{->}[0,1]^(.575){\pmKum{\localfield{\cyc}}{p^\infty}}&\tors{\HH{1}{\localfield{\cyc}}{\tors{\EC}{p^\infty}}}{p}\ar@{->}[0,1]\ar@{->}[0,1]^-{\alpha}&\tors{\Bigl(\HH{1}{\localfield{\cyc}}{\tors{\EC}{p^\infty}}/\image\pmKum{\localfield{\cyc}}{p^\infty}\Bigr)}{p}
}\end{split}
\end{equation}
The first vertical arrow is an isomorphism in light of~\eqref{eq:EC/p=EC[p]}, and the second vertical arrow is an isomorphism thanks to~\ref{point:prop:H[p]=H(p)_infinite_global}. The snake lemma implies that $\pmfakeKum[\ssprimeid{\empty}]{\cyc}$ is injective and $\coker(\pmfakeKum[\ssprimeid{\empty}]{\cyc})=\coker(\alpha)$. To show that $\alpha$ is surjective observe that the bottom row in~\eqref{diag:Kumm_modp_ss} is the beginning of the $\operatorname{Tor}^i_{\ZZ_p}(-,\ZZ/p)$-sequence of the tautological exact sequence
\[
0\longrightarrow\pmEC(\localfield{\cyc})\otimes\QQ_p/\ZZ_p\xrightarrow{\pmKum{\localfield{\cyc}}{p^\infty}}
\HH{1}{\localfield{\cyc}}{\tors{\EC}{p^\infty}}\longrightarrow\HH{1}{\localfield{\cyc}}{\tors{\EC}{p^\infty}}/\image\pmKum{\localfield{\cyc}}{p^\infty}\longrightarrow 0
\]
and therefore $\coker(\alpha)$ is contained in
\[
\Bigr(\pmEC(\localfield{\cyc})\otimes\QQ_p/\ZZ_p\Bigl)\otimes\ZZ/p=\Bigr(\pmEC(\localfield{\cyc})\otimes\QQ_p/\ZZ_p\Bigl)/p\Bigr(\pmEC(\localfield{\cyc})\otimes\QQ_p/\ZZ_p\Bigl).
\]
Since $\pmEC(\localfield{\cyc})\otimes\QQ_p/\ZZ_p$ is divisible, the above module is trivial, establishing the surjectivity of $\alpha$ and thus of $\pmfakeKum[\ssprimeid{\empty}]{\cyc}$. This finishes the proof of the proposition.
\end{proof}

Let us now pass to Selmer groups. We refer to \cite[Chapter~2]{CoaSuj10} for generalities on Iwasawa theory for elliptic curves over cyclotomic extensions and, in particular, for the definitions of the groups $\fullSel{\basecyc}{\tors{\EC}{p^\infty}}$ in the ordinary case.
\begin{definition}\label{def:multi_signed_Sel}The multi-signed residual Selmer group $\pmSel*{\tors{\EC}{p}}{\basecyc}$ is defined as
\[
\pmSel{\tors{\EC}{p}}{\basecyc}=\varinjlim_{\res{}{}}\pmSel{\tors{\EC}{p}}{\basefield[n]}
\]
and the usual multi-signed Selmer group is defined as
\[
\KobSel{\tors{\EC}{p^\infty}}{\basecyc}=\varinjlim_{\res{}{}}\KobSel{\tors{\EC}{p^\infty}}{\basefield[n]}.
\]
\end{definition}
The groups $\pmSel{\tors{\EC}{p}}{\basecyc}$ are discrete $\resLambda$-modules, whose Pontryagin duals are compact, finitely generated over $\resLambda$. Similarly, the groups $\KobSel{\tors{\EC}{p^\infty}}{\basecyc}$ are discrete, cofinitely generated $\IwLambda$-modules. For $\genprimeid\in\primeset$, denote by $\modJJ{\genprimeid}$ the $\resLambda$-module
\[
\modJJ{\genprimeid}=\begin{cases}
\displaystyle{
\bigoplus_{\GenPrimeId\mid\badprimeid{\empty}}\HH{1}{\genlocal[\GenPrimeId]{\cyc}}{\tors{\EC}{p}}}&\text{ if }\genprimeid=\badprimeid{\empty}\in\badprimeset\\[2em]
\displaystyle{\bigoplus_{\GenPrimeId\mid\ordprimeid{\empty}}\HH{1}{\genlocal[\GenPrimeId]{\cyc}}{\tors{\modEC}{p}}}&\text{ if }\genprimeid=\ordprimeid{\empty}\in\ordprimeset\\[2em]
\HH{1}{\ssloccyc[i]}{\tors{\EC}{p}}/\image(\pmKum{\ssloccyc[i]}{p})&\text{ if }\genprimeid=\ssprimeid{i}\in\ssprimeset.
\end{cases}\]
\begin{remark} It is not \emph{a priori} obvious that the groups $\modJJ{\genprimeid}$ depend only upon the isomorphism class of $\tors{\EC}{p}$, as the notation would suggest. Indeed, for a supersingular prime $\ssprimeid{}\in\ssprimeset$, the definition of the map $\pmKum{\ssloccyc[i]}{p}$ is given in terms of the full torsion module $\tors{\EC}{p^\infty}$. Below, we will see that in fact the groups $\modJJ{\genprimeid}$ only depend upon $\tors{\EC}{p}$.
\end{remark}
Similarly, define $\JJ{\genprimeid}$ as the $\IwLambda$-module
\[
\JJ{\genprimeid}=\bigoplus_{\GenPrimeId\mid\genprimeid}\HH{1}{\genlocal[\GenPrimeId]{\cyc}}{\tors{\EC}{p^\infty}}/\image(\pmKum{\genlocal[\GenPrimeId]{\cyc}}{p^\infty})
\]
where, as in Proposition~\ref{prop:H[p]=H(p)_infinite}, we set $\pmKum{\genlocal[\GenPrimeId]{\cyc}}{p^\infty}=\fullKum{\genlocal[\GenPrimeId]{\cyc}}{p^\infty}$ for all $\GenPrimeId\mid\genprimeid\in\primeset\setminus\ssprimeset$. As in Definition~\ref{def:pm-Sel}, when a vector $\multisign\in\{+,-\}^{\ordssprimeset}$ is fixed, denote simply by $\modJJ{\genprimeid}$ (\resp by $\JJ{\genprimeid}$) the module $\modJJ[+]{\genprimeid}$ or $\modJJ[-]{\genprimeid}$ (\resp the module $\JJ[+]{\genprimeid}$ or $\JJ[-]{\genprimeid}$) according to the $i$th component $\multisign_i\in\{+,-\}$. Then the multi-signed residual Selmer group and the usual multi-signed Selmer group sit in the following exact sequences:
\begin{align}\label{eq:def_modJJ}
0\longrightarrow\pmSel{\tors{\EC}{p}}{\basecyc}\longrightarrow\HH{1}{\GS[\cyc]}{\tors{\EC}{p}}\xrightarrow{\pmlambda{\cyc}{p}}\displaystyle{\bigoplus_{\genprimeid\in\primeset}\modJJ{\genprimeid}}\\
\intertext{\resp}
\label{eq:def_fullJJ}
0\longrightarrow\KobSel{\tors{\EC}{p^\infty}}{\basecyc}\longrightarrow\HH{1}{\GS[\cyc]}{\tors{\EC}{p^\infty}}\xrightarrow{\pmlambda{\cyc}{p^\infty}}\displaystyle{\bigoplus_{\genprimeid\in\primeset}\bigl(\JJ{\genprimeid}\bigr)}.
\end{align}
As a consequence of the results in Section~\ref{section:signed}, we can now check that the local conditions $\modJJ{\genprimeid}$ depend only on the residual representation $\tors{\EC}{p}$. This is immediate at all $\badprimeid{}\in\badprimeset$, and follows from Proposition~\ref{prop:raynaud}-\ref{item:prop:raynaud_modE=modE} at all primes $\ordprimeid{}\in\ordprimeset$ by taking inductive limit along the cyclotomic tower. Concerning primes $\ssprimeid{}\in\ssprimeset$, this independence follows from combining Proposition~\ref{prop:universal_sequence} with the isomorphism in Proposition~\ref{prop:H[p]=H(p)_infinite}-\ref{point:prop:H[p]=H(p)_infinite_ss}. The fact that the local conditions depend only on the residual representation implies that the same holds for the multi-signed residual Selmer group. We record this as
\begin{corollary}\label{cor:Sel_indep} Let $\EC_1,\EC_2$ be two elliptic curves over $\basefield$ satisfying~\ref{hyp_1} such that $\tors{(\EC_1)}{p}\cong\tors{(\EC_2)}{p}$. Then, the sets $\ssprimeset[1]$ and $\ssprimeset[2]$ of primes of supersingular reduction for $\EC_1$ and $\EC_2$ coincide and for every vector of signs $\multisign\in\{+,-\}^{\ordssprimeset}$,
\[
\pmSel{\tors{(\EC_1)}{p}}{\basecyc}\cong\pmSel{\tors{(\EC_2)}{p}}{\basecyc}.
\]
\end{corollary}
\begin{proof}
The equality $\ssprimeset[1]=\ssprimeset[2]$ is a consequence of Proposition~\ref{prop:raynaud}-\ref{item:prop:raynaud_Sss=Sord}, and the isomorphism of the corresponding multi-signed residual Selmer groups follows from the above discussion.
\end{proof}

Kim considers in \cite{Kim09} the primitive and non-primitive Selmer groups along the lines of \cite{GreVat00}. Corollary~\ref{cor:H[p]=H(p)_infinite_selmer} below is the analogue of~\cite[Proposition~2.10]{Kim09} for the multi-signed residual Selmer groups. In order to state it, let us introduce a final notation. For all $\GenPrimeId\mid\genprimeid\in\badprimeset\cup\ordprimeset$,  let $\ordcycprimes{v}$ be the number of primes $\GenPrimeId$ lying above $\genprimeid$ in $\basecyc$. Further, for $v=\badprimeid{\empty}\in\badprimeset$ and $\BadPrimeId{}\mid\badprimeid{}$ in $\basecyc$, denote by $\cyclocal{\BadPrimeId{}}$ the first layer of the cyclotomic $\ZZ_p$-extension $\genlocal[\cyc]{\BadPrimeId{}}/\genlocal{\badprimeid{}}$. Note that $\cyclocal{\BadPrimeId{}}$ is also the unique unramified extension of degree $p$ of~$\genlocal{\badprimeid{}}$.

Recall that, for any place $\genprimeid$, the residue field at that place is denoted by $\finitefield[\genprimeid]$.
\begin{definition}\label{def:defect} For all $\badprimeid{\empty}\in\badprimeset$, choose a place $\BadPrimeId{\empty}$ of $\basecyc$ above $\badprimeid{\empty}$. We define the defect of $\EC$ as
\begin{equation*}
\defect{\EC}:=\sum_{\badprimeid{\empty}\in\badprimeset}\ordcycprimes{\badprimeid{\empty}}\cdot\dim_{\finitefield}\tors{\EC(\cyclocal{\BadPrimeId{\empty}})}{p}+\sum_{\ordprimeid{\empty}\in\ordprimeset}\ordcycprimes{\ordprimeid{\empty}}\dim_{\finitefield}\tors{\modEC(\finitefield[\ordprimeid{\empty}])}{p}\leq 2\sum_{\badprimeid{\empty}\in\badprimeset}\ordcycprimes{\badprimeid{\empty}}+\sum_{\ordprimeid{\empty}\in\ordprimeset}\ordcycprimes{\ordprimeid{\empty}}.
\end{equation*}
The fact that $\dim_{\finitefield}\tors{\EC(\cyclocal{\BadPrimeId{\empty}})}{p}$ is independent of $\BadPrimeId{\empty}\mid\badprimeid{\empty}$ follows from $\basecyc/\basefield$ being Galois, since $\EC$ is defined over $\basefield$.
\end{definition}
\begin{corollary}\label{cor:H[p]=H(p)_infinite_selmer} 
There are injections
\[
\varphi^\multisign\colon\pmSel{\tors{\EC}{p}}{\basecyc}\longhookrightarrow\tors{\KobSel{\tors{\EC}{p^\infty}}{\basecyc}}{p}
\]
whose cokernel is finite, of dimension $\dim_{\finitefield}\coker(\varphi^\multisign)\leq\defect{\EC}$. In particular,
\[
\corank[\resLambda]\pmSel{\tors{\EC}{p}}{\basecyc}=\corank[\resLambda]\tors{\KobSel{\tors{\EC}{p^\infty}}{\basecyc}}{p}.
\]
Moreover, when $\pmlambda{\cyc}{p}$ is surjective, $\dim_{\finitefield}\coker(\varphi^\multisign)=\defect{\EC}$, independently of the vector $\multisign$.
\end{corollary}
\begin{proof}
Consider the commutative diagram
\begin{equation}\label{diag:Selp_KobSel}\begin{split}\xymatrix{
0\ar@{->}[0,1]&\pmSel{\tors{\EC}{p}}{\basecyc}\ar@{->}[0,1]\ar@{->}[1,0]^{\varphi^\multisign}&\HH{1}{\GS[\cyc]}{\tors{\EC}{p}}\ar@{->}[0,1]^-{\pmlambda{\cyc}{p}}\ar@{->}[1,0]^{\rotatebox{270}{$\cong$}}&\displaystyle{\bigoplus_{\genprimeid\in\primeset}\modJJ{\genprimeid}}\ar@{->}[1,0]^{\bigoplus{\varphi_{\genprimeid}^\pm}}\\
0\ar@{->}[0,1]&\tors{\KobSel{\tors{\EC}{p^\infty}}{\basecyc}}{p}\ar@{->}[0,1]&\tors{\HH{1}{\GS[\cyc]}{\tors{\EC}{p^\infty}}}{p}\ar@{->}[0,1]^-{\pmlambda{\cyc}{p^\infty}}&\displaystyle{\bigoplus_{\genprimeid\in\primeset}\tors{\bigl(\JJ{\genprimeid}\bigr)}{p}}
}\end{split}\end{equation}
The central vertical arrow is an isomorphism thanks to Proposition~\ref{prop:H[p]=H(p)_infinite}-\ref{point:prop:H[p]=H(p)_infinite_global}. The local arrows $\varphi_{\genprimeid}^\pm$ can be decomposed as
\[
\varphi_{\genprimeid}^\pm=\bigoplus_{\GenPrimeId\mid\genprimeid}\varphi_{\GenPrimeId}^\pm
\]
and each $\varphi_{\GenPrimeId}^\pm$ is induced by the corresponding arrow $\fakeKum[\GenPrimeId]{\cyc}$ of Proposition~\ref{prop:H[p]=H(p)_infinite}. More precisely,
\[
\varphi_{\GenPrimeId}^\pm=\begin{cases}
\fakeKum[\GenPrimeId]{\cyc}\colon\HH{1}{\genlocal[\GenPrimeId]{\cyc}}{\tors{\EC}{p}}\longrightarrow \tors{\HH{1}{\genlocal[\GenPrimeId]{\cyc}}{\tors{\EC}{p^\infty}}}{p}=\tors{\Bigl(\HH{1}{\genlocal[\GenPrimeId]{\cyc}}{\tors{\EC}{p^\infty}}/\image(\pmKum{\genlocal[\GenPrimeId]{\cyc}}{p^\infty})\Bigr)}{p}&\text{ if}\quad\genprimeid=\badprimeid{\empty}\in\badprimeset\\[2em]
\modfakeKum[\GenPrimeId]{\cyc}\colon\HH{1}{\genlocal[\GenPrimeId]{\cyc}}{\tors{\modEC}{p}}\longrightarrow\tors{\HH{1}{\genlocal[\GenPrimeId]{\cyc}}{\tors{\modEC}{p^\infty}}}{p}= \tors{\Bigl(\HH{1}{\genlocal[\GenPrimeId]{\cyc}}{\tors{\EC}{p^\infty}}/\image(\pmKum{\genlocal[\GenPrimeId]{\cyc}}{p^\infty})\Bigr)}{p}&\text{ if}\quad\genprimeid=\ordprimeid{\empty}\in\ordprimeset\\[2em]
\pmfakeKum[\ssprimeid]{\cyc}\colon\HH{1}{\ssloccyc}{\tors{\EC}{p}}/\image(\pmKum{\ssloccyc}{p})\overset{\cong}{\longrightarrow} \tors{\Bigl(\HH{1}{\ssloccyc}{\tors{\EC}{p^\infty}}/\image(\pmKum{\ssloccyc}{p^\infty})\Bigr)}{p}
&\text{ if}\quad\genprimeid=\ssprimeid{\empty}\in\ssprimeset
\end{cases}
\]
Indeed, the first equality
\[
\tors{\HH{1}{\genlocal[\GenPrimeId]{\cyc}}{\tors{\EC}{p^\infty}}}{p}=\tors{\Bigl(\HH{1}{\genlocal[\GenPrimeId]{\cyc}}{\tors{\EC}{p^\infty}}/\image(\pmKum{\genlocal[\GenPrimeId]{\cyc}}{p^\infty})\Bigr)}{p}
\]
follows from the fact that $\image(\pmKum{\genlocal[\GenPrimeId]{\cyc}}{p^\infty})=0$, since $\EC(\genlocal[\GenPrimeId]{\cyc})\otimes\QQ_p/\ZZ_p=0$ when $\genprimeid\notin\pprimeset$. The second equality
\[
\tors{\HH{1}{\genlocal[\GenPrimeId]{\cyc}}{\tors{\modEC}{p^\infty}}}{p}= \tors{\Bigl(\HH{1}{\genlocal[\GenPrimeId]{\cyc}}{\tors{\EC}{p^\infty}}/\image(\pmKum{\genlocal[\GenPrimeId]{\cyc}}{p^\infty})\Bigr)}{p}
\]
follows from~\eqref{eq:CoaGre:modp}. Similarly, the fact that $\pmfakeKum[\ssprimeid{\empty}]{\cyc}$ takes values in $\tors{\Bigl(\HH{1}{\ssloccyc}{\tors{\EC}{p^\infty}}/\image(\pmKum{\ssloccyc}{p^\infty})\Bigr)}{p}$ and is an isomorphism follows from Proposition~\ref{prop:H[p]=H(p)_infinite}-\ref{point:prop:H[p]=H(p)_infinite_ss}.

By applying the snake lemma to~\eqref{diag:Selp_KobSel} we see that
\begin{equation}\label{eq:coker_Selp_Sbad:Sord}
\coker(\varphi^\multisign)\subseteq\bigoplus_{\BadPrimeId{\empty}\mid\badprimeid{\empty}\in\badprimeset}\ker\bigl(\fakeKum[\BadPrimeId{}]{\cyc}\bigr)\oplus\bigoplus_{\GenPrimeId\mid\ordprimeid{\empty}\in\ordprimeset}\ker\bigl(\modfakeKum[\GenPrimeId]{\cyc}\bigr)
\end{equation}
and the inclusion in~\eqref{eq:coker_Selp_Sbad:Sord} is an equality when $\pmlambda{\cyc}{p}$ is surjective. 

It follows from~\eqref{eq:coker_Selp_Sbad:Sord} and Proposition~\ref{prop:H[p]=H(p)_infinite}-\ref{point:prop:H[p]=H(p)_infinite_bad}--\ref{point:prop:H[p]=H(p)_infinite_ord} that
\begin{equation*}
\dim_{\finitefield}\coker(\varphi^\multisign)\leq\dim_{\finitefield}\bigoplus_{\BadPrimeId{}\mid\badprimeid{\empty}\in\badprimeset}\dim_{\finitefield}\tors{\EC(\genlocal[\cyc]{\GenPrimeId})}{p}+\dim_{\finitefield}\bigoplus_{\GenPrimeId\mid\ordprimeid{\empty}\in\ordprimeset}\dim_{\finitefield}\tors{\modEC(\finitefield[\GenPrimeId])}{p}
\end{equation*}
which is an equality if $\pmlambda{\cyc}{p}$ is surjective. To prove the statement of the corollary, we need to show that
\begin{align}\label{eq:cor:H[p]=H(p)_bad}
&\dim_{\finitefield}\tors{\EC(\genlocal[\cyc]{\BadPrimeId{}})}{p}=\tors{\EC(\cyclocal{\BadPrimeId{\empty}})}{p}\qquad&&\text { if }\BadPrimeId{}\mid\badprimeid{}\in\badprimeset
\intertext{and}\label{eq:cor:H[p]=H(p)_ord}
&\dim_{\finitefield}\tors{\modEC(\finitefield[\GenPrimeId])}{p}=\dim_{\finitefield}\tors{\modEC(\finitefield[\ordprimeid{\empty}])}{p}\qquad&&\text { for all }\GenPrimeId\mid\ordprimeid{}\in\ordprimeset.
\end{align}
The equality in~\eqref{eq:cor:H[p]=H(p)_ord} simply follows from the fact that the inertia degree of every $p$-adic prime is $1$ along the $\ZZ_p$-cyclotomic $\basecyc/\basefield$, whence $\finitefield[\GenPrimeId]=\finitefield[\ordprimeid{}]$.

Concerning~\eqref{eq:cor:H[p]=H(p)_bad}, we argue as follows. If $\dim_{\finitefield}\tors{\EC(\genlocal{\badprimeid{}})}{p}=2$, then the full torsion of $\EC(\overline{\genlocal{\badprimeid{}}})$ is  already defined over $\genlocal{\badprimeid{}}$, hence $\tors{\EC(\genlocal{\badprimeid{}})}{p}=\tors{\EC(\cyclocal{\BadPrimeId{}})}{p}=\tors{\EC(\genlocal[\cyc]{\BadPrimeId{}})}{p}$. If $\dim_{\finitefield}\tors{\EC(\genlocal{\badprimeid{}})}{p}=0$, then the $p$-group $\tors{\EC(\genlocal[\cyc]{\BadPrimeId{}})}{p}$ has no non-zero fixed point under the action of the pro-$p$-group $\Gamma=\Gal{\genlocal[\cyc]{\BadPrimeId{}}}{\genlocal{\badprimeid{}}}$, and must be trivial. Hence $\tors{\EC(\genlocal[\cyc]{\BadPrimeId{}})}{p}=0=\tors{\EC(\cyclocal{\BadPrimeId{}})}{p}$. Finally, if $\dim_{\finitefield}\tors{\EC(\genlocal{\badprimeid{}})}{p}=1$ but $\dim_{\finitefield}\tors{\EC(\genlocal[\cyc]{\BadPrimeId{}})}{p}=2$, we need to show that $\dim_{\finitefield}\tors{\EC(\cyclocal{\BadPrimeId{}})}{p}=2$. The assumption that $\dim_{\finitefield}\tors{\EC(\genlocal[\cyc]{\BadPrimeId{}})}{p}=2$ implies that the action of $\Gamma$ in $\operatorname{Aut}\bigl(\tors{\EC(\genlocal[\cyc]{\GenPrimeId})}{p}\bigr)$ induces, upon fixing a basis, a $2$-dimensional linear representation
\[
\varrho\colon\Gamma\longrightarrow \operatorname{GL}_2(\finitefield).
\] 
As $\operatorname{GL}_2(\finitefield)$ contains no element of order $p^2$, $\varrho$ factors through $\Gamma/\Gamma^p=\Gal{\cyclocal{\BadPrimeId{}}}{\genlocal{\badprimeid{}}}$, so $\Gamma^p$ acts trivially on $\tors{\EC(\genlocal[\cyc]{\BadPrimeId{}})}{p}$. This implies that $\tors{\EC(\genlocal[\cyc]{\BadPrimeId{}})}{p}=\tors{\EC(\cyclocal{\BadPrimeId{}})}{p}$ also in this case, and concludes the proof of the corollary.
\end{proof}

\subsection{Cassels--Poiutou--Tate exact sequence and Iwasawa cohomology}\label{subsec:poitou-tate} For $n\in \NN\cup\{\cyc\}$, let $\dualSel{\tors{\EC}{p^\infty}}{\basefield[n]}$ denote the Pontryagin duals
\[
\dualSel{\tors{\EC}{p^\infty}}{\basefield[n]}=\Hom_{\ZZ_p}\bigl(\fullSel{\tors{\EC}{p^\infty}}{\basefield[n]},\QQ_p/\ZZ_p\bigr).
\]
These modules clearly admit multi-signed versions, defined as
\[
\KobdualSel{\tors{\EC}{p^\infty}}{\basefield[n]}=\Hom_{\ZZ_p}(\KobSel{\tors{\EC}{p^\infty}}{\basefield[n]},\QQ_p/\ZZ_p).
\]
Similarly, the duals of the multi-signed residual Selmer groups are defined as
\[
\dualpmSel{\tors{\EC}{p}}{\basefield[n]}=\Hom_{\ZZ_p}(\pmSel{\tors{\EC}{p}}{\basefield[n]},\QQ_p/\ZZ_p).
\]
Both $\dualSel{\tors{\EC}{p^\infty}}{\basecyc}$ and $\KobdualSel{\tors{\EC}{p^\infty}}{\basecyc}$ are finitely generated compact $\IwLambda$-modules and it follows from Corollary~\ref{cor:H[p]=H(p)_infinite_selmer} that the $\resLambda$-modules $\dualpmSel{\tors{\EC}{p}}{\basefield[n]}$ are finitely generated. Further, Corollary~\ref{cor:Sel_indep} implies that they only depend upon the isomorphism class of $\tors{\EC}{p}$. As a last piece of notation, suppose that $\field\in\{\localfield{n},\genlocal{\genprimeid},\basefield[n]\}$ for some $0\leq n<\infty$, and retain notation from~\eqref{not:S_and_primes}: in particular, $\widetilde{\field}=\overline{\QQ_p}$ if $\field=\localfield{n}$, $\widetilde{\field}=\overline{\genlocal{\genprimeid}}$ if $\field=\genlocal{\genprimeid}$ and $\widetilde{\field}=\baseS$ if $\field=\basefield$. Let $M$ be a compact $\ZZ_p$-module with a continuous $\Gal{\widetilde{\field}}{\field}$-action. The \emph{Iwasawa cohomology} modules $\HHiw{i}{\field}{M}$ (for all $i\geq 1$) are defined as the inverse limit, with respect to corestriction maps
\[
\HHiw{i}{\field}{M}=\varprojlim_{\field\subseteq \field'\subseteq\fieldcyc}\HH{i}{\fieldS/\field'}{M}.
\]
The reader is referred to \cite[\S3.1]{Per92} for generalities about Iwasawa cohomology. In particular, the above $\IwLambda$-modules are known to be trivial for $i\neq 1,2$.

A fundamental tool for the study of the Iwasawa theory of Selmer groups is the Cassels--Poitou--Tate exact sequence, for which we refer to \cite[Theorem~1.5]{CoaSuj10} and which we now briefly recall. Fix~$n\in\NN$ and consider the self-dual module $M=\tors{\EC}{p}$. Put
\[
W_{\genprimeid}=\begin{cases}
0&\text{ if }\genprimeid=\badprimeid{\empty}\in\badprimeset\\
\image(\pmKum{\localfield[i]{n}}{p})&\text{ if }\genprimeid=\ssprimeid{i}\in\ssprimeset\quad (1\leq i\leq \ordssprimeset)\\
\HH{1}{\localfield[\ordprimeid{i}]{n}}{\tors{(\ordformalEC_{\ordprimeid{i}})}{p}}&\text { if }\genprimeid=\ordprimeid{i}\in\ordprimeset \quad (1\leq i\leq\ordordprimeset)
\end{cases}\]
where $\ordformalEC_{\ordprimeid{}}$ denotes the formal group of $\EC/\rint{\genlocal{\ordprimeid{}}}$, and the sign $\pm$ depends on the $i$-th component of $\multisign$. These coincide with the local conditions in Definition~\ref{def:pm-Sel}, so that the group $\pmSel{\tors{\EC}{p}}{\basefield[n]}$ sits in the exact sequence
\[
0\longrightarrow\pmSel{\tors{\EC}{p}}{\basefield[n]}\longrightarrow\HH{1}{\GS[n]}{\tors{\EC}{p}}\longrightarrow\bigoplus_{\genprimeid\in\primeset}\HH{1}{\genlocal[\genprimeid]{n}}{\tors{\EC}{p}}/W_{\genprimeid}.
\]
For all $\genprimeid\in\primeset$, let $W_{\genprimeid}^\perp$ denote the orthogonal complement of $W_{\genprimeid}$ in the Tate pairing
\[
\HH{1}{\genlocal[\genprimeid]{n}}{\tors{\EC}{p}}\times \HH{1}{\genlocal[\genprimeid]{n}}{\tors{\EC}{p}}\longrightarrow\QQ/\ZZ
\]
and define $\CPTSel{\tors{\EC}{p}}{\basefield}$ as the kernel
\[
0\longrightarrow\CPTSel{\tors{\EC}{p}}{\basefield[n]}\longrightarrow\HH{1}{\GS[n]}{\tors{\EC}{p}}\longrightarrow\bigoplus_{\genprimeid\in\primeset}\HH{1}{\genlocal[\genprimeid]{n}}{\tors{\EC}{p}}/W_{\genprimeid}^\perp.
\]
This gives the Cassels--Poitou--Tate exact sequence
\begin{equation}\begin{aligned}\label{seq:pmcasselspoitoutate_finite}
0\longrightarrow\pmSel{\tors{\EC}{p}}{\basefield[n]}\longrightarrow \HH{1}{\GS[n]}{\tors{\EC}{p}}\longrightarrow \bigoplus_{\genprimeid\in\primeset}&\HH{1}{\genlocal[\genprimeid]{n}}{\tors{\EC}{p}}/W_{\genprimeid}
{\longrightarrow}\pontr{\Bigl(\CPTSel{\tors{\EC}{p}}{\basefield[n]}\Bigr)}\longrightarrow\\
&\longrightarrow \HH{2}{\GS[n]}{\tors{\EC}{p}}\longrightarrow \bigoplus_{\GenPrimeId\mid\genprimeid\in \primeset}\HH{2}{\genlocal[n]{\genprimeid}}{\tors{\EC}{p}}\longrightarrow 0
\end{aligned}\end{equation}
where the final~$0$ comes from Proposition~\ref{prop:no_p_torsion}. Similarly, put
\[
U_{\genprimeid}=\begin{cases}
0&\text{ if }\genprimeid=\badprimeid{\empty}\in\badprimeset\\
\image(\pmKum{\localfield[i]{n}}{p^\infty})&\text{ if }\genprimeid=\ssprimeid{i}\in\ssprimeset\quad (1\leq i\leq \ordssprimeset)\\
\HH{1}{\localfield[\ordprimeid{i}]{n}}{\tors{(\ordformalEC_{\ordprimeid{i}})}{p^\infty}}&\text { if }\genprimeid=\ordprimeid{i}\in\ordprimeset \quad (1\leq i\leq\ordordprimeset)
\end{cases}\]
and define $U_{\genprimeid}^\perp\subseteq\HH{1}{\genlocal[n]{\genprimeid}}{\tate{\EC}}$ to be the orthogonal complement of $U_{\genprimeid}$. Defining $\KobCPT{\tate{\EC}}{\basefield}$ as the kernel
\[
0\longrightarrow\KobCPT{\tate{\EC}}{\basefield[n]}\longrightarrow\HH{1}{\GS[n]}{\tate{\EC}}\longrightarrow\bigoplus_{\genprimeid\in\primeset}\HH{1}{\genlocal[\genprimeid]{n}}{\tate{\EC}}/U_{\genprimeid}^\perp.
\]
we obtain the Cassels--Poitou--Tate exact sequence
\begin{equation}\begin{aligned}\label{seq:pmcasselspoitoutate_finite_full}
0\longrightarrow\KobSel{\tors{\EC}{p^\infty}}{\basefield[n]}\longrightarrow \HH{1}{\GS[n]}{\tors{\EC}{p^\infty}}\longrightarrow \bigoplus_{\genprimeid\in\primeset}&\HH{1}{\genlocal[\genprimeid]{n}}{\tors{\EC}{p^\infty}}/U_{\genprimeid}
{\longrightarrow}\pontr{\Bigl(\KobCPT{\tate{\EC}}{\basefield[n]}\Bigr)}\longrightarrow\\
&\longrightarrow \HH{2}{\GS[n]}{\tors{\EC}{p^\infty}}\longrightarrow \bigoplus_{\GenPrimeId\mid\genprimeid\in \primeset}\HH{2}{\genlocal[n]{\genprimeid}}{\tors{\EC}{p^\infty}}\longrightarrow 0
\end{aligned}\end{equation}

In order to study to the limit, as $n\to\infty$ along the cyclotomic tower, of the Cassels--Poitou--Tate sequences, consider the groups
\begin{equation*}
\CPTprojlim{\tors{\EC}{p}}{\basecyc}=\varprojlim\CPTSel{\tors{\EC}{p}}{\basefield[n]}\subseteq \HHiw{1}{\basefield}{\tors{\EC}{p}}.
\end{equation*}
as well as
\begin{equation*}
\KobprojSel{\tate{\EC}}{\basecyc}=\varprojlim\KobCPT{\tate{\EC}}{\basefield[n]}\subseteq \HHiw{1}{\basefield}{\tate{\EC}}.
\end{equation*}
where inverse limits are taken with respect to corestriction maps. These are, respectively, a $\resLambda$-module and a $\IwLambda$-module and their relevance for our study comes from the following observation (\cfr \cite[Lemma~2.6 and Remark~2.7]{LimSuj18b}, where the case of an elliptic curve with ordinary reduction at $p$ is considered):
\begin{lemma}\label{lemma:duality_lim}
There are isomorphisms
\[
\pontr{\CPTprojlim{\tors{\EC}{p}}{\basecyc}}\cong\varinjlim\pontr{\Bigl(\CPTSel{\tors{\EC}{p}}{\basefield[n]}\Bigr)}\qquad\text{and}\qquad \pontr{\KobprojSel{\tate{\EC}}{\basecyc}}\cong\varinjlim\pontr{\Bigl(\KobCPT{\tate{\EC}}{\basefield[n]}\Bigr)}
\]
where the direct limits are taken with respect to the dual of corestriction maps. Moreover, the group $\CPTprojlim{\tors{\EC}{p}}{\basecyc}$ is free as $\resLambda$-module and $\KobprojSel{\tate{\EC}}{\basecyc}$ has no non-zero torsion $\IwLambda$-submodules.
\end{lemma}
\begin{proof} The first isomorphism simply follows from the definition, since
\begin{align*}
\pontr{\CPTprojlim{\tors{\EC}{p}}{\basecyc}}&=\Hom\Bigl(\varprojlim\CPTSel*{\tors{\EC}{p}}{\basefield[n]},\QQ_p/\ZZ_p\Bigr)\\
&=\varinjlim\Hom\Bigl(\CPTSel*{\tors{\EC}{p}}{\basefield[n]},\QQ_p/\ZZ_p\Bigr)\\
&=\varinjlim\pontr{\Bigl(\CPTSel{\tors{\EC}{p}}{\basefield[n]}\Bigr)}.
\end{align*}
The second isomorphism is analogous.

To show that $\CPTprojlim{\tors{\EC}{p}}{\basecyc}$ is a free $\resLambda$-module, note that Jannsen's spectral sequence~\cite[Corollary~13]{Jan14} takes the form
\[
E_{2}^{p,q}=\Ext_{\resLambda}^p\Bigl(\Hom\bigl(\HH{q}{\GS[\cyc]}{\tors{\EC}{p}},\finitefield\bigr),\resLambda\Bigr)\Longrightarrow \HHiw{p+q}{\basefield}{\tors{\EC}{p}}.
\]
By Proposition~\ref{prop:no_p_torsion}, $E_2^{p,0}=0$ for all $p\geq 0$. Therefore we have  $E_{\infty}^{1,0} = E_2^{1,0}=0$ and $E_2^{0,1} =E_{\infty}^{0,1},$ and it follows that
\[
E_2^{0,1}=\Hom_{\resLambda}\bigl(\Hom_{\finitefield}(\HH{1}{\GS[\cyc]}{\tors{\EC}{p}},\finitefield),\resLambda\bigr)\cong \HHiw{1}{\basefield}{\tors{\EC}{p}}.
\]
In particular, the first Iwasawa cohomology group of $\tors{\EC}{p}$ is torsion-free over $\resLambda$, and hence free since~$\resLambda$ is a PID. By taking inverse limit, with respect to corestriction, of the inclusions $\CPTSel{\tors{\EC}{p}}{\basefield[n]}\hookrightarrow \HH{1}{\GS[n]}{\tors{\EC}{p}}$, we obtain an injection
\[
\CPTprojlim{\tors{\EC}{p}}{\basecyc}\xhookrightarrow{\phantom{\longrightarrow}}\HHiw{1}{\basefield}{\tors{\EC}{p}}.
\]
Since $\resLambda$ is a principal ideal domain, a submodule of a free module is itself free, thereby proving the claim concerning $\CPTprojlim{\tors{\EC}{p}}{\basecyc}$. Passing to the Galois representation $\tors{\EC}{p^\infty}$, Jannsen's spectral sequence~\cite[Theorem~1]{Jan14} is
\[
E_{2}^{p,q}=\Ext_{\IwLambda}^p\Bigl(\Hom\bigl(\HH{q}{\GS[\cyc]}{\tors{\EC}{p^\infty}},\QQ_p/\ZZ_p\Bigr),\IwLambda\Bigr)\Longrightarrow \HHiw{p+q}{\basefield}{\tate{\EC}}.
\]
Again, Proposition~\ref{prop:no_p_torsion}, yields $E_2^{p,0}=0$ for all $p\geq 1$ and the same argument as above shows that the first Iwasawa cohomology group of $\tate{\EC}$ is torsion-free over $\IwLambda$. In particular, the same holds for its submodule $\KobprojSel{\tate{\EC}}{\basecyc}$.
\end{proof}
In~\S~\ref{subsubsection:final_notation}, there is a summary of the notation for the various Selmer groups, and their signed versions, for both $\tors{\EC}{p}$ and for $\tors{\EC}{p^\infty}$.

Given two subextensions $\basefield\subseteq\basefield[n]\subseteq \basefield[m]\subseteq \basecyc$, consider the corresponding exact sequences~\eqref{seq:pmcasselspoitoutate_finite}. The restriction map on cohomology induces morphisms between the first three (\resp the last two) terms. A standard argument in local Tate duality shows that connecting the fourth terms via the Pontryagin dual of corestriction
\[
\pontr{(\cores{}{})}\colon \pontr{\Bigl(\CPTSel{\tors{\EC}{p}}{\basefield[n]}\Bigr)}\longrightarrow\pontr{\Bigl(\CPTSel{\tors{\EC}{p}}{\basefield[m]}\Bigr)}\qquad(m\geq n\geq 0)
\]
of~\eqref{seq:pmcasselspoitoutate_finite}, gives commutative diagrams of exact sequences. By taking the direct limit over $n$, Lemma~\ref{lemma:duality_lim} gives exact sequences
\begin{equation}\begin{aligned}\label{seq:pmcpt_cyc}
0\longrightarrow\pmSel{\tors{\EC}{p}}{\basecyc}\longrightarrow \HH{1}{\GS[\cyc]}{\tors{\EC}{p}}\xrightarrow{\pmlambda{\basecyc}{p}} \bigoplus_{\genprimeid\in \primeset}&\modJJ{\genprimeid}
\longrightarrow\pontr{\CPTprojlim{\tors{\EC}{p}}{\basecyc}}\longrightarrow\\
\longrightarrow& \HH{2}{\GS[\cyc]}{\tors{\EC}{p}}\longrightarrow \bigoplus_{\GenPrimeId\mid\genprimeid\in \primeset}\HH{2}{\basecyc[\GenPrimeId]}{\tors{\EC}{p}}\longrightarrow 0
\end{aligned}\end{equation}
where the morphism $\pmlambda{\basecyc}{p}$ and the $\resLambda$-modules $\modJJ{\genprimeid}$ were introduced in~\eqref{eq:def_modJJ}. Replacing the Galois representation $\tors{\EC}{p}$ by $\tors{\EC}{p^\infty}$ we obtain the analogous exact sequence
\begin{equation}\begin{aligned}\label{seq:pmcpt_cyc_infty}
0\longrightarrow\KobSel{\tors{\EC}{p^\infty}}{\basecyc}\longrightarrow \HH{1}{\GS[\cyc]}{\tors{\EC}{p^\infty}}\xrightarrow{\pmlambda{\basecyc}{p^\infty}} \bigoplus_{\genprimeid\in \primeset}&\JJ{\genprimeid}
\longrightarrow\pontr{\KobprojSel{\tate{\EC}}{\basecyc}}\longrightarrow\\
\longrightarrow& \HH{2}{\GS[\cyc]}{\tors{\EC}{p^\infty}}\longrightarrow \bigoplus_{\GenPrimeId\mid\genprimeid\in \primeset}\HH{2}{\basecyc[\GenPrimeId]}{\tors{\EC}{p^\infty}}\longrightarrow 0.
\end{aligned}\end{equation}

We go back to the study of $\pmSel{\tors{\EC}{p}}{\basecyc}$. Unlike the ordinary case, the full Selmer group
\begin{equation}\label{eq:def_Selmer_infty}
\begin{split}
 \fullSel{\tors{\EC}{p^\infty}}{\basecyc}&=\ker\Bigl(\HH{1}{\GS[\infty]}{\tors{\EC}{p^\infty}}\longrightarrow \bigoplus_{\GenPrimeId\mid\genprimeid\in \primeset}\tors{\HH{1}{\basecyc[\GenPrimeId]}{\EC}}{p^\infty}\Bigr)\\
&=\ker\biggl(\HH{1}{\GS[\infty]}{\tors{\EC}{p^\infty}}\longrightarrow \bigoplus_{\GenPrimeId\mid\genprimeid\in \primeset\setminus\ssprimeset}\tors{\HH{1}{\basecyc[\GenPrimeId]}{\EC}}{p^\infty}\oplus\bigoplus_{i=1}^{\ordssprimeset}\HH{1}{\ssloccyc[i]}{\tors{\EC}{p^\infty}}/\image(\fullKum{\ssloccyc[i]}{p^\infty})\biggr)
\end{split}\end{equation}
is not $\IwLambda$-cotorsion, in general. Indeed, as is discussed in the proof of \cite[Theorem~2.6]{CoaSuj10}, each local term $\tors{\HH{1}{\ssloccyc[i]}{\EC}}{p^\infty}$ has $\IwLambda$-corank equal to $0$ or $1$ depending on the reduction type of $\modEC/\finitefield$ and this implies that $\fullSel{\tors{\EC}{p^\infty}}{\basecyc}$ is not $\IwLambda$-cotorsion in the supersingular case. In the ordinary reduction case, Mazur asked in \cite[\S6]{Maz72} whether $\fullSel{\tors{\EC}{p^\infty}}{\basecyc}$ is co-torsion. This is known to be true if $\fullSel{\tors{\EC}{p^\infty}}{\basefield}$ is finite or if $\basefield=\QQ$ (see \cite[Theorem~2.8 and Theorem~2.18]{CoaSuj10}). 

In the supersingular setting, both Perrin-Riou and Kobayashi reduce the size of the kernels in~\eqref{eq:def_Selmer_infty}  by replacing $\image(\fullKum{\ssloccyc}{p^\infty})$ with the smaller subgroup $\image(\pmKum{\ssloccyc}{p^\infty})$, at supersingular primes. This has the effect that the corresponding \emph{signed} Selmer group is potentially a cotorsion module over the Iwasawa algebra. We follow the same strategy, replacing $\tors{\EC}{p^\infty}$ by $\tors{\EC}{p}$, and replacing signed Selmer groups by their \emph{fine multi-signed} residual versions.

\subsubsection{Notation}\label{subsubsection:final_notation}
We gather here the notation concerning Selmer groups introduced above, together with a reference pointing at its definition. 
\begin{itemize}
\item The local conditions $\JJ{\genprimeid}$ and $\modJJ{\genprimeid}$ are the defined right after Definition~\ref{def:multi_signed_Sel} and are subgroups of global cohomology groups of $\tors{\EC}{p^\infty}$ and of $\tors{\EC}{p}$, respectively.
\item These local conditions yield the definition of the (signed) Selmer groups $\KobSel{\tors{\EC}{p^\infty}}{\basecyc}$ and  $\pmSel{\tors{\EC}{p}}{\basecyc}$. The rationale beneath the choice of the letter \emph{R} for the Selmer group of the representation $\tors{\EC}{p}$ is to hint at a \emph{residual} Selmer group.
\item Taking Pontryagin duals of the above groups, one obtains the groups 
\[
\KobdualSel{\tors{\EC}{p^\infty}}{\basefield[n]}=\Hom_{\ZZ_p}(\KobSel{\tors{\EC}{p^\infty}}{\basefield[n]},\QQ_p/\ZZ_p).
\]
and
\[
\dualpmSel{\tors{\EC}{p}}{\basefield[n]}=\Hom_{\ZZ_p}(\pmSel{\tors{\EC}{p}}{\basefield[n]},\QQ_p/\ZZ_p).
\]
defined at the beginning of~\S~\ref{subsec:poitou-tate}.
\item By local Tate duality, the local conditions $\JJ{\genprimeid}$ and $\modJJ{\genprimeid}$ give rise to dual local conditions. These, in turn, define dual compact Selmer groups $\KobCPT{\tate{\EC}}{\basefield[n]}$ and $\CPTSel{\tors{\EC}{p}}{\basefield[n]}$, appearing in equations~\eqref{seq:pmcasselspoitoutate_finite_full} and~\eqref{seq:pmcasselspoitoutate_finite}, respectively. Again, the choice of letters \emph{S} and \emph{R} reflect the fact that these are the compact versions of the ``full'' Selmer group and of the residual Selmer group, respectively.
\item Taking inverse limits, with respect to corestriction maps, over the cyclotomic tower of the compact groups $\KobCPT{\tate{\EC}}{\basefield[n]}$ and $\CPTSel{\tors{\EC}{p}}{\basefield[n]}$ defines the $\IwLambda$-modules and $\KobprojSel{\tate{\EC}}{\basecyc}$ $\CPTprojlim{\tors{\EC}{p}}{\basecyc}$ (the latter is in fact a $\resLambda$-module). Once more, the choice of letters is coherent with the previous rationale. The relation between Pontryagin and Tate duality is summarised in Lemma~\ref{lemma:duality_lim}
\end{itemize}
\subsection{Rank computation}To approach $\pmSel{\tors{\EC}{p}}{\basecyc}$, the main objects of study will be the maps $\pmlambda{\basecyc}{p}$ and $\pmlambda{\basecyc}{p^\infty}$ appearing in the exact sequences~\eqref{seq:pmcpt_cyc} and~\eqref{seq:pmcpt_cyc_infty}. We will ultimately relate the surjectivity of $\pmlambda{\basecyc}{p}$ to the structure of $\KobdualSel{\tors{\EC}{p^\infty}}{\basecyc}$ as a $\IwLambda$-module (see Theorem~\ref{thm:main_result}). In this direction, the following hypothesis (which is the multi-sign analogue of condition (vi) in~\cite[Theorem~1.3]{KitOts18}) is crucial: 
\begin{equation}\tag*{\textup{\textbf{Hyp~2}}$_\multisign$}\label{hyp_2}
\begin{parbox}{0.85\textwidth}{
The $\IwLambda$-module $\KobdualSel{\tors{\EC}{p^\infty}}{\basecyc}$ 
is torsion and hence admits two structural Iwasawa invariants, which we denote by $\iwlambda$ and $\iwmu$.
}\end{parbox}\end{equation}

The exact sequence which is crucial in our approach is~\eqref{seq:pmcpt_cyc}.
Concerning the term $\HH{2}{\GS[\cyc]}{\tors{\EC}{p}}$ in it, Coates and the second author have proposed in \cite{CoaSuj05} the following
\begin{labelledconj}{Conjecture~\textup{A}}{conj:A}
$\HH{2}{\GS[\cyc]}{\tors{\EC}{p}}=0$.
\end{labelledconj}
The original formulation of~\ref{conj:A} \ibid~is that the dual fine Selmer group $\dualfineSel{\EC}{\basecyc}$ over $\basecyc$ (see~\cite[\S3]{CoaSuj05} for its definition) is a finitely generated $\ZZ_p$-module. Note that this is equivalent to $\dualfineSel{\EC}{\basecyc}$ being $\IwLambda$-torsion, and having $\iwmu[\empty]$-invariant equal to $0$. The next proposition relates the two formulations, and shows that~\ref{conj:A} implies the Weak Leopoldt Conjecture (see Remark~\ref{rmk:rmks_propConj} below).
\begin{proposition}\label{prop:conjectures} \ref{conj:A} implies that $\HH{2}{\GS[\cyc]}{\tors{\EC}{p^\infty}}=0$. Moreover, if $\EC$ satisfies~\ref{hyp_2}, and if $\iwmu=0$,
then~\ref{conj:A} holds.
\end{proposition}
\begin{proof}
Taking $\GS[\cyc]$-cohomology of the exact sequence~\eqref{eq:false_kummer} yields
\[
\HH{2}{\GS[\cyc]}{\tors{\EC}{p}}\longrightarrow\HH{2}{\GS[\cyc]}{\tors{\EC}{p^\infty}}\overset{\cdot p}{\longrightarrow}\HH{2}{\GS[\cyc]}{\tors{\EC}{p^\infty}}\longrightarrow 0
\]
where the surjection comes from the fact that $\Gal{\overline{\basefield}}{\basecyc}$ has cohomological dimension $2$ (see~\cite[Theorem~10.11.3 and~Proposition~3.3.5]{NeuSchWin08}). Therefore, if $\HH{2}{\GS[\cyc]}{\tors{\EC}{p}}=0$ then multiplication by $p$ is injective on $\HH{2}{\GS[\cyc]}{\tors{\EC}{p^\infty}}$. On the other hand, every class in $\HH{2}{\GS[\cyc]}{\tors{\EC}{p^\infty}}$ has finite $p$-power order, hence multiplication by $p$ is injective if and only if $\HH{2}{\GS[\cyc]}{\tors{\EC}{p^\infty}}=0$.

Suppose now that $\KobdualSel{\tors{\EC}{p^\infty}}{\basecyc}$ is $\IwLambda$-torsion and $\iwmu=0$. The dual fine Selmer group $\dualfineSel{\EC}{\basecyc}$ is defined in~\cite[(42)]{CoaSuj05} as the Pontryagin dual of
\[
\ker\Bigl(\HH{1}{\GS[\cyc]}{\tors{\EC}{p^\infty}}\longrightarrow\bigoplus_{\GenPrimeId\mid\genprimeid\in \primeset}\HH{1}{\basecyc[\GenPrimeId]}{\tors{\EC}{p^\infty}}\Bigr).
\]
Since this kernel injects into $\KobSel{\tors{\EC}{p^\infty}}{\basecyc}$ by~\eqref{eq:def_fullJJ}, we obtain a surjection
\[
\KobdualSel{\tors{\EC}{p^\infty}}{\basecyc}\longtwoheadrightarrow \dualfineSel{\EC}{\basecyc}.
\]
Our assumptions imply then that $\dualfineSel{\EC}{\basecyc}$ is a torsion $\IwLambda$-module with trivial $\iwmu[\empty]$-invariant, which is the formulation of~\cite[Conjecture~A]{CoaSuj05}. We are thus left to show that if $\dualfineSel{\EC}{\basecyc}$ is $\IwLambda$-torsion and has trivial $\iwmu[\empty]$-invariant, then $\HH{2}{\GS[\cyc]}{\tors{\EC}{p}}=0$. But Greenberg shows in~\cite[Proposition~4.1.6]{Gre11} that the vanishing of $\HH{2}{\GS[\cyc]}{\tors{\EC}{p}}$ is equivalent to the $p$-torsion subgroup $\tors{\dualfineSel{\EC}{\basecyc}}{p}$ being finite, and this is certainly the case when $\dualfineSel{\EC}{\basecyc}$ is $\IwLambda$-torsion and has trivial $\iwmu[\empty]$-invariant.
\end{proof}
\begin{remark}\label{rmk:rmks_propConj}\leavevmode
\begin{enumerate}[label=\arabic*.]
\item \label{rmk:rmks_propConjA:wlc}The vanishing of $\HH{2}{\GS[\cyc]}{\tors{\EC}{p^\infty}}$ is known as the \emph{Weak Leopoldt Conjecture} (see \cite[p.~348]{Sch85} and \cite[Conjecture~3]{Gre89}). If $\basefield=\QQ$, it holds by \cite[Theorem~12.4]{Kat04}, at least for $\primeset=\{p\}$; the case for general $\primeset=\badprimeset\cup\{p\}$ can be deduced from Kato's result by combining the exact sequence of~\cite[p.~33, (1.4.3)]{Per95} with Jannsen's spectral sequence from~\cite[Theorem~1]{Jan14}. Over an arbitrary base $\basefield$, if $\fullSel{\EC}{\basefield}$ is finite, the vanishing can be proven by combining~\cite[Proposition~1.9]{CoaSuj10} with the Hochschild--Serre spectral sequence.
\item \label{rmk:rmks_propConjA_modular} It is clear that the validity of~\ref{conj:A} depends only upon the isomorphism class of the Galois representation $\tors{\EC}{p}$.
\end{enumerate}
\end{remark}

Let us now pass to the study of $\resLambda$-coranks of some cohomology groups, which will turn out to be a key step in the proof of our main result. In Lemma~\ref{lemma:divisible_stuff}, global cohomology and local cohomology at primes where $\EC$ does not have supersingular reduction are considered. Then, in Lemma~\ref{lemma:plus/minus_coatesgreenberg}, supersingular primes are treated.
\begin{lemma}\label{lemma:divisible_stuff}
Let $\genprimeid\in\primeset\setminus\ssprimeset$ and let $\GenPrimeId\mid \genprimeid$ be a place in $\basecyc$ that lies above $\genprimeid$. Then,
\begin{enumerate}[label=\roman*)]
\item \label{point:lemma:divisible_stuff_bad} If $\genprimeid=\badprimeid{\empty}\in\badprimeset$, the Pontryagin duals $\HH{1}{\genlocal[\cyc]{\GenPrimeId}}{\tors\EC{p^\infty}}$ have $\iwmu[]$ invariant equal to $0$.
\item \label{point:lemma:divisible_stuff_ord} If $\genprimeid=\ordprimeid{\empty}\in\ordprimeset$, the Pontryagin duals of $\HH{1}{\genlocal[\cyc]{\GenPrimeId}}{\tors{\EC}{p^\infty}}/\image(\fullKum{\genlocal[\cyc]{\GenPrimeId}}{p^\infty})$ have $\iwmu[]$ invariant equal to $0$.
\item \label{point:lemma:divisible_stuff_global} Assuming~\ref{conj:A}, the Pontryagin dual of $\HH{1}{\GS[\cyc]}{\tors{\EC}{p^\infty}}$ has trivial $\iwmu[]$ invariant as well.
\end{enumerate}
As a consequence,
\begin{align}\label{eq:lemma:divisible_bad}
\corank[\resLambda]\HH{1}{\genlocal[\cyc]{\GenPrimeId}}{\tors{\EC}{p}}
&=\corank\HH{1}{\genlocal[\cyc]{\GenPrimeId}}{\tors{\EC}{p^\infty}}
,&\quad\GenPrimeId\mid\badprimeid{\empty}\in\badprimeset\\
\label{eq:lemma:divisible_ord}
\corank[\resLambda]\HH{1}{\genlocal[\cyc]{\GenPrimeId}}{\tors{\modEC}{p}}
&=\corank\Bigl(\HH{1}{\genlocal[\cyc]{\GenPrimeId}}{\tors{\EC}{p^\infty}}\big/\image(\fullKum{\genlocal[\cyc]{\GenPrimeId}}{p^\infty})\Bigr)&\,\GenPrimeId\mid\ordprimeid{\empty}\in\ordprimeset\\
\intertext{and, assuming~\ref{conj:A},}
\label{eq:lemma:divisible_global}
\corank[\resLambda]\HH{1}{\GS[\cyc]}{\tors{\EC}{p}}&=\corank\HH{1}{\GS[\cyc]}{\tors{\EC}{p^\infty}}.
\end{align}
\end{lemma}
\begin{proof} We start with local cohomology, and let $\genprimeid\in\primeset\setminus\ssprimeset$ be any prime. Greenberg proves in \cite[Propositions~1 and~2]{Gre89} that the groups $\HH{1}{\genlocal[\cyc]{\GenPrimeId}}{\tors{\EC}{p^\infty}}$ are cofinitely generated: this implies, in particular, that their quotients $\HH{1}{\genlocal[\cyc]{\GenPrimeId}}{\tors{\EC}{p^\infty}}/\image(\fullKum{\genlocal[\cyc]{\GenPrimeId}}{p^\infty})$ are cofinitely generated as well. Moreover, we claim that the exact sequence~\eqref{eq:false_kummer} induces
\[
\HH{1}{\genlocal[\cyc]{\GenPrimeId}}{\tors{\EC}{p^\infty}}\overset{\cdot p}{\longrightarrow}\HH{1}{\genlocal[\cyc]{\GenPrimeId}}{\tors{\EC}{p^\infty}}\longrightarrow \HH{2}{\genlocal[\cyc]{\GenPrimeId}}{\tors{\EC}{p}}=0.
\]
The $\operatorname{H}^2$-term in the above sequence vanishes because $\Glocal[{\genlocal[\cyc]}]{\GenPrimeId}$ has $p$-cohomological dimension~$1$, as observed in the proof of Proposition~\ref{prop:H[p]=H(p)_infinite}.

The fact that multiplication by $p$ is surjective on $\HH{1}{\genlocal[\cyc]{\GenPrimeId}}{\tors{\EC}{p^\infty}}$ shows that this module is $p$-divisible, and thus the same holds for the quotient $\HH{1}{\genlocal[\cyc]{\GenPrimeId}}{\tors{\EC}{p^\infty}}/\image(\fullKum{\genlocal[\cyc]{\GenPrimeId}}{p^\infty})$. Observe now that this divisibility is equivalent to their Pontryagin duals having no $p$-torsion and, in particular, to having trivial $\iwmu[]$ invariant. This establishes points~\ref{point:lemma:divisible_stuff_bad} and~\ref{point:lemma:divisible_stuff_ord}.

When~\ref{conj:A} holds, the same argument as above shows that multiplication by $p$ is surjective on $\HH{1}{\GS[\cyc]}{\tors{\EC}{p^\infty}}$, whence~\ref{point:lemma:divisible_stuff_global}.

By Proposition~\ref{prop:H[p]=H(p)_infinite}-\ref{point:prop:H[p]=H(p)_infinite_global} (\resp~Proposition~\ref{prop:H[p]=H(p)_infinite}-\ref{point:prop:H[p]=H(p)_infinite_bad}), the Pontryagin duals of the $\resLambda$-modules $\HH{1}{\GS[\cyc]}{\tors{\EC}{p}}$ and $\tors{\HH{1}{\GS[\cyc]}{\tors{\EC}{p^\infty}}}{p}$  
(\resp $\HH{1}{\genlocal[\cyc]{\GenPrimeId}}{\tors{\EC}{p}}$ and  $\tors{\HH{1}{\genlocal[\cyc]{\GenPrimeId}}{\tors{\EC}{p^\infty}}}{p}$ for some $\GenPrimeId\mid\badprimeid{\empty}\in\badprimeset$) have the same rank. Now equations~\eqref{eq:lemma:divisible_bad} and~\eqref{eq:lemma:divisible_global} follow from assertions~\ref{point:lemma:divisible_stuff_bad} and~\ref{point:lemma:divisible_stuff_global}, respectively, along with the structure theorem for finitely generated $\IwLambda$-modules. Similarly, combining Proposition~\ref{prop:H[p]=H(p)_infinite}-\ref{point:prop:H[p]=H(p)_infinite_ord} with~\ref{point:lemma:divisible_stuff_ord} yields~\eqref{eq:lemma:divisible_ord}.
\end{proof}
We finish the study of $\resLambda$-coranks of cohomology groups by analysing what happens at supersingular primes. Our argument is the analogue, modulo $p$, of~\cite[Proposition~3.32]{KitOts18}.
\begin{lemma}\label{lemma:plus/minus_coatesgreenberg} For each choice of sign $\pm$, the $\resLambda$-module
\[
\pontr{\Bigl(\HH{1}{\ssloccyc}{\tors{\EC}{p}}/\image(\pmKum{\ssloccyc}{p})\Bigr)}
\]
is finitely generated and free of rank $1$.
\end{lemma}
\begin{proof} The statement will follow once we prove that
\begin{equation}\label{eq:quotient_free}
\pontr{\Bigl(\HH{1}{\ssloccyc}{\tors{\EC}{p^\infty}}/\image(\pmKum{\ssloccyc}{p^\infty})\Bigr)}\cong\IwLambda,
\end{equation}
thanks to Proposition~\ref{prop:H[p]=H(p)_infinite}-\ref{point:prop:H[p]=H(p)_infinite_ss}.

The freeness claimed in~\eqref{eq:quotient_free} follows from \cite[Lemma~3.31]{KitOts18}. Indeed,
\[
\pontr{\Bigl(\HH{1}{\ssloccyc}{\tors{\EC}{p^\infty}}/\image(\pmKum{\ssloccyc}{p^\infty})\Bigr)}=\ker\Bigl(
\pontr{\HH{1}{\ssloccyc}{\tors{\EC}{p^\infty}}}
\longtwoheadrightarrow
\pontr{\image(\pmKum{\ssloccyc}{p^\infty})}\Bigr)
\]
where $\pontr{\HH{1}{\ssloccyc}{\tors{\EC}{p^\infty}}}$ is $\IwLambda$-free of rank $2$, as proven in~\cite[Corollary~1]{Gre89}, and $\pontr{\image(\pmKum{\ssloccyc}{p^\infty})}$ is of $\IwLambda$-rank equal to $1$ and has no non-trivial finite $\IwLambda$-submodules, as follows from~\cite[Proposition~3.28~(for $\chi=1$)]{KitOts18}.
\end{proof}
The following Proposition is essentially well-known in the ordinary case, and it has already been proven by Iovita--Pollack in the supersingular case under the assumption that $\EC$ is defined over $\QQ$, and $p$ splits completely in $\basefield/\QQ$ (see {\cite[Proposition~6.1]{IovPol06}}).
\begin{proposition}\label{prop:ranks}
Suppose that~\ref{conj:A} holds for $\EC/\basefield$. Then we have
\begin{align}\label{eq:ranks_IwLambda}
\corank[\IwLambda]\HH{1}{\GS[\cyc]}{\tors{\EC}{p^\infty}}=&\sum_{\ordprimeid{\empty}\in\ordprimeset}\corank[\IwLambda]\JJ{\ordprimeid{}}+\sum_{i=1}^{\ordssprimeset}\corank[\IwLambda]\JJ{\ssprimeid{i}}\\
\label{eq:ranks_ResLambda}
\corank[\resLambda]\HH{1}{\GS[\cyc]}{\tors{\EC}{p}}=&\sum_{\ordprimeid{\empty}\in\ordprimeset}\corank[\resLambda]\modJJ{\ordprimeid{}}+\sum_{i=1}^{\ordssprimeset}\corank[\resLambda]\modJJ{\ssprimeid{i}}.
\end{align}
Moreover,
\[
\sum_{\badprimeid{\empty}\in\badprimeset}\corank[\IwLambda]\JJ{\badprimeid{}}=\sum_{\badprimeid{\empty}\in\badprimeset}\corank[\resLambda]\modJJ{\badprimeid{}}=0.
\]
\end{proposition}
\begin{proof} The proof is an adaptation of~\cite[Proof of Theorem~2.6]{CoaSuj10}. We first compute the left-hand sides of both~\eqref{eq:ranks_IwLambda} and~\eqref{eq:ranks_ResLambda}. In \cite[Proposition~3]{Gre89} Greenberg proves that both $\HH{1}{\GS[\cyc]}{\tors{\EC}{p^\infty}}$ and $\HH{2}{\GS[\cyc]}{\tors{\EC}{p^\infty}}$ are co-finitely generated over $\IwLambda$ and further 
\begin{equation*}
\corank\HH{1}{\GS[\cyc]}{\tors{\EC}{p^\infty}}-\corank\HH{2}{\GS[\cyc]}{\tors{\EC}{p^\infty}}=2r_2+\sum_{\genprimeid\text{ real place}}d_{\genprimeid}^-.
\end{equation*}
Here $r_2$ is the number of complex places of $\basefield$ and, for each real place $\genprimeid$ of $\basefield$, we denote by $d_{\genprimeid}^-$ the dimension of the $(-1)$-eigenspace for a complex conjugation above $\genprimeid$ acting on $\tate{\EC}\otimes\QQ_p$. By Proposition~\ref{prop:conjectures}, the $\operatorname{H}^2$-term vanishes and, by the Galois invariance of the Weil pairing, we know that $d_{\genprimeid}^-=1$ for all real $\genprimeid$. Hence,
\begin{equation*}
\corank\HH{1}{\GS[\cyc]}{\tors{\EC}{p^\infty}}=[\basefield:\QQ]=\degbase.
\end{equation*}
Now~\eqref{eq:lemma:divisible_global} of Lemma~\ref{lemma:divisible_stuff} implies
\begin{equation*}
\corank[\resLambda]\HH{1}{\GS[\cyc]}{\tors{\EC}{p}}=\corank\HH{1}{\GS[\cyc]}{\tors{\EC}{p^\infty}}=\degbase.
\end{equation*}

Passing to the computation of the local coranks, let first $\ordprimeid{\empty}\in\ordprimeset$. By~\cite[\S 2.13]{CoaSuj10} (which applies here, thanks to our convention that $\pmKum{\basecyc[\GenPrimeId]}{p^\infty}=\fullKum{\basecyc[\GenPrimeId]}{p^\infty}$ when $\ordprimeid{\empty}\in\ordprimeset$, together with~\eqref{eq:CoaGre:modp}) 
we know
\begin{equation}\label{eq:greenberg_ordp_infty}
\corank\bigoplus_{\GenPrimeId\mid\ordprimeid{\empty}}\Bigl(\HH{1}{\genlocal[\cyc]{\GenPrimeId}}{\tors{\EC}{p^\infty}}/\image(\pmKum{\basecyc[\GenPrimeId]}{p^\infty})\Bigr)=[\genlocal{\ordprimeid{\empty}}:\QQ_p].
\end{equation}
Hence equation~\eqref{eq:lemma:divisible_ord} of Lemma~\ref{lemma:divisible_stuff} yields
\begin{equation}\label{eq:greenberg_ordp}
\corank[\resLambda]\bigoplus_{\GenPrimeId\mid\ordprimeid{\empty}}\HH{1}{\genlocal[\cyc]{\GenPrimeId}}{\tors{\modEC}{p}}
=\corank[\resLambda]\modJJ{\ordprimeid{}}=[\genlocal{\ordprimeid{\empty}}:\QQ_p].
\end{equation}
Consider now a prime~$\ssprimeid{i}\in\ssprimeset$. Lemma~\ref{lemma:plus/minus_coatesgreenberg} implies that
\begin{equation}\label{eq:greenberg_ssp}
\corank[\resLambda]\bigl(\HH{1}{\ssloccyc[i]}{\tors{\EC}{p}}/\image(\pmKum{\ssloccyc}{p})\bigr)=\corank[\resLambda]\modJJ{\ssprimeid{i}}=1.
\end{equation}
Combining~\eqref{eq:greenberg_ordp_infty} and~\eqref{eq:quotient_free}, we find
\begin{equation*}
\sum_{\ordprimeid{\empty}\in\ordprimeset}\corank[\IwLambda]\JJ{\ordprimeid{}}+\sum_{i=1}^{\ordssprimeset}\corank[\IwLambda]\JJ{\ssprimeid{i}}=\degbase
\end{equation*}
Similarly,~\eqref{eq:greenberg_ordp} and~\eqref{eq:greenberg_ssp} together imply 
\begin{equation*}
\sum_{\ordprimeid{\empty}\in\ordprimeset}\corank[\resLambda]\modJJ{\ordprimeid{}}+\sum_{i=1}^{\ordssprimeset}\corank[\resLambda]\modJJ{\ssprimeid{i}}=\degbase
\end{equation*}
and this establishes equations~\eqref{eq:ranks_IwLambda} and~\eqref{eq:ranks_ResLambda}

Finally, suppose that $\badprimeid{\empty}$ is a prime in $\badprimeset$ and let $\BadPrimeId{}$ be an extension of $\badprimeid{\empty}$ to $\basecyc$. Greenberg proves in \cite[Proposition~2]{Gre89} that the $\IwLambda$-module $\HH{1}{\genlocal[\cyc]{\BadPrimeId{}}}{\tors{\EC}{p^\infty}}$ is cotorsion. Hence
\begin{equation*}
\corank[\IwLambda]\JJ{\badprimeid{}}=\sum_{\BadPrimeId{}\mid\badprimeid{}}\corank\HH{1}{\genlocal[\cyc]{\BadPrimeId{}}}{\tors{\EC}{p^\infty}}
=0
\end{equation*}
and~\eqref{eq:lemma:divisible_bad} of Lemma~\ref{lemma:divisible_stuff} implies $\sum_{\BadPrimeId{}\mid\badprimeid{}}\corank[\resLambda]\HH{1}{\genlocal[\cyc]{\BadPrimeId{}}}{\tors{\EC}{p}}=\sum_{\BadPrimeId{}\mid\badprimeid{}}\corank\HH{1}{\genlocal[\cyc]{\BadPrimeId{}}}{\tors{\EC}{p^\infty}}
=0$ as well. This completes the proof of the proposition.
\end{proof}
\subsection{Main results} We are now in a position to state and prove our main result. Recall the exact sequence
\begin{equation}\tag{\ref{seq:pmcpt_cyc}}\begin{aligned}
0\longrightarrow\pmSel{\tors{\EC}{p}}{\basecyc}\longrightarrow \HH{1}{\GS[\cyc]}{\tors{\EC}{p}}\xrightarrow{\pmlambda{\basecyc}{p}} \bigoplus_{\genprimeid\in \primeset}&\modJJ{\genprimeid}
\longrightarrow\pontr{\CPTprojlim{\tors{\EC}{p}}{\basecyc}}\longrightarrow\\
\longrightarrow& \HH{2}{\GS[\cyc]}{\tors{\EC}{p}}\longrightarrow \bigoplus_{\GenPrimeId\mid\genprimeid\in \primeset}\HH{2}{\basecyc[\GenPrimeId]}{\tors{\EC}{p}}\longrightarrow 0
\end{aligned}\end{equation}
and consider the projection
\[
\ppr\colon\bigoplus_{\genprimeid\in \primeset}\modJJ{\genprimeid}\rightarrow\bigoplus_{\genprimeid\in \pprimeset}\modJJ{\genprimeid}.
\]
Define $\plocallambda{\tors{\EC}{p}}$, or simply $\plocallambda{p}$, as the composition $\ppr\kern-.35em\circ\kern.25em\pmlambda{\basecyc}{p}$:
\begin{theorem}\label{thm:main_result} Under our standing assumption~\ref{hyp_1}, suppose that~\ref{hyp_2} holds as well. Then the following assertions are equivalent:
\begin{enumerate}[label=\alph*\textup{)}]
\item\label{thm:main_result:surj} $\pmlambda{\basecyc}{p}$ is surjective and~\ref{conj:A} holds;
\item\label{thm:main_result:loc_surj} $\plocallambda{p}$ is surjective and~\ref{conj:A} holds;
\item\label{thm:main_result:cotorsion} $\pmSel{\tors{\EC}{p}}{\basecyc}$ is $\resLambda$-cotorsion;
\item\label{thm:main_result:mu} $\KobdualSel{\tors{\EC}{p^\infty}}{\basecyc}$ has trivial $\iwmu[\empty]$-invariant.
\end{enumerate}
\end{theorem}
\begin{proof}To show that~\ref{thm:main_result:surj}~$\Rightarrow$~\ref{thm:main_result:cotorsion}, take Pontryagin duals of the short exact sequence
\begin{equation*}
0\longrightarrow\pmSel{\tors{\EC}{p}}{\basecyc}\longrightarrow \HH{1}{\GS[\cyc]}{\tors{\EC}{p}}\xrightarrow{\pmlambda{\basecyc}{p}} \bigoplus_{\genprimeid\in \primeset}\modJJ{\genprimeid}
\longrightarrow 0
\end{equation*}
to obtain
\begin{equation}\label{eq:surj_xi:shortES}
0\longrightarrow \pontr{\bigoplus_{\genprimeid\in \primeset}\modJJ{\genprimeid}}\xrightarrow{\phantom{\kern.5em}\pontr{\pmlambda{\basecyc}{p}}\phantom{\kern.5em}} \pontr{\HH{1}{\GS[\cyc]}{\tors{\EC}{p}}}\longrightarrow\dualpmSel{\tors{\EC}{p}}{\basecyc}
\longrightarrow 0.
\end{equation}
By Proposition~\ref{prop:ranks}, the first two terms have the same $\resLambda$-rank, so the third is $\resLambda$-torsion and~\ref{thm:main_result:cotorsion} follows. Also, when~\ref{thm:main_result:surj} holds, the composition $\plocallambda{p}=\ppr\kern-.35em\circ\kern.25em\pmlambda{\basecyc}{p}$ is surjective, yielding~\ref{thm:main_result:surj}~$\Rightarrow$~\ref{thm:main_result:loc_surj}.

To show that~\ref{thm:main_result:cotorsion}~and~\ref{thm:main_result:mu} are equivalent, we first observe that a finitely generated torsion $\IwLambda$-module $M$ has trivial $\iwmu[\empty]$-invariant if and only if $M/pM$ is a torsion $\resLambda$-module. On the other hand, taking Pontryagin duals of the injection of Corollary~\ref{cor:H[p]=H(p)_infinite_selmer} shows that the kernel of
\[
\pontr{\tors{\Bigr(\KobSel{\tors{\EC}{p^\infty}}{\basecyc}\Bigl)}{p}}=\KobdualSel{\tors{\EC}{p^\infty}}{\basecyc}/p	\KobdualSel{\tors{\EC}{p^\infty}}{\basecyc}\longtwoheadrightarrow\dualpmSel{\tors{\EC}{p}}{\basecyc}
\]
is finite, showing the equivalence between~\ref{thm:main_result:cotorsion}~and~\ref{thm:main_result:mu}. Therefore,~\ref{thm:main_result:loc_surj}~$\Leftarrow$~\ref{thm:main_result:surj}~$\Rightarrow$~\ref{thm:main_result:cotorsion}~$\Leftrightarrow$~\ref{thm:main_result:mu}.

We are left with the implications~\ref{thm:main_result:cotorsion}~$\Rightarrow$~\ref{thm:main_result:surj} and~\ref{thm:main_result:loc_surj}~$\Rightarrow$~\ref{thm:main_result:surj}. Since~\ref{thm:main_result:cotorsion}~$\Rightarrow$~\ref{thm:main_result:mu}~$\Rightarrow$~\ref{conj:A} by Proposition~\ref{prop:conjectures}, and~\ref{thm:main_result:loc_surj} contains~\ref{conj:A}, we can assume from now on that $\HH{2}{\GS[\cyc]}{\tors{\EC}{p}}=0$. In particular, the sequence~\eqref{seq:pmcpt_cyc} becomes
\begin{equation}\label{seq:pmcpt_cyc_ConjA}
0\longrightarrow\pmSel{\tors{\EC}{p}}{\basecyc}\longrightarrow \HH{1}{\GS[\cyc]}{\tors{\EC}{p}}\xrightarrow{\pmlambda{\basecyc}{p}} \bigoplus_{\genprimeid\in \primeset}\modJJ{\genprimeid}
\longrightarrow\pontr{\CPTprojlim{\tors{\EC}{p}}{\basecyc}}=\coker(\pmlambda{\basecyc}{p})\longrightarrow 0
\end{equation}
and Proposition~\ref{prop:ranks} yields
\begin{equation*}
\corank[\resLambda]\bigr(\pmSel{\tors{\EC}{p}}{\basecyc}\bigl)=\corank[\resLambda]\bigl(\CPTprojlim{\tors{\EC}{p}}{\basecyc}\bigr).
\end{equation*}
Assuming~\ref{thm:main_result:cotorsion}, it follows that $\CPTprojlim{\tors{\EC}{p}}{\basecyc}$ is $\resLambda$-torsion. On the other $\CPTprojlim{\tors{\EC}{p}}{\basecyc}$ is $\resLambda$-free, in light of Lemma~\ref{lemma:duality_lim}, and to be $\resLambda$-torsion it must be trivial, establishing \ref{thm:main_result:cotorsion}~$\Rightarrow$~\ref{thm:main_result:surj}.

Finally, assume that $\plocallambda{p}$ is surjective and consider the commutative triangle
\[\xymatrix{
\HH{1}{\GS[\cyc]}{\tors{\EC}{p}}\ar@{->}[0,2]^-{\pmlambda{\basecyc}{p}}\ar@{->}[1,2]_{\plocallambda{p}}&&\displaystyle{\bigoplus_{\badprimeid{}\in\badprimeset}\modJJ{\badprimeid{}}\oplus\bigoplus_{\genprimeid\in\pprimeset}\modJJ{\genprimeid}}\ar@{->}[1,0]^-{\ppr}\\
&&
\displaystyle{\bigoplus_{\genprimeid\in\pprimeset}\modJJ{\genprimeid{}}}.
}\]
It induces an exact sequence
\begin{equation*}
\ker (\ppr)=\bigoplus_{\badprimeid{}\in\badprimeset}\modJJ{\badprimeid{}}\longrightarrow\coker(\pmlambda{\basecyc}{p})\longrightarrow\coker(\plocallambda{p})=0
\end{equation*}
which implies that $\coker(\pmlambda{\basecyc}{p})$ is cotorsion, thanks to Proposition~\ref{prop:ranks}.
Since $\coker(\pmlambda{\basecyc}{p})$ is isomorphic to $\pontr{\CPTprojlim{\tors{\EC}{p}}{\basecyc}}$ by~\eqref{seq:pmcpt_cyc_ConjA}, and $\CPTprojlim{\tors{\EC}{p}}{\basecyc}$ is free by Lemma~\ref{lemma:duality_lim}, this forces $\coker(\pmlambda{\basecyc}{p})=0$, establishing the final implication~\ref{thm:main_result:loc_surj}~$\Rightarrow$~\ref{thm:main_result:surj}.
\end{proof}
\begin{remark} Note that~\ref{conj:A} is pivotal to the proof and plays the role of the Weak Leopoldt Conjecture for the residual representation~$\tors{\EC}{p}$.
\end{remark}
As an application of Theorem~\ref{thm:main_result}, we obtain a result along the lines of Greenberg--Vatsal's work \cite[Theorem~1.4]{GreVat00} in the supersingular setting. Results somewhat similar to Theorem~\ref{thm:mu=mu} below, again in the supersingular setting, have been obtained by Kim in~\cite[Corollary~2.13]{Kim09} and by Hatley--Lei in~\cite[Theorem~4.6]{HatLei19}, by different methods.

\begin{theorem}\label{thm:mu=mu}
Let $\EC_1,\EC_2$ be two elliptic curves defined over $\basefield$, satisfying hypothesis~\ref{hyp_1} and such that the residual Galois representations $\tors{(\EC_1)}{p}$ and $\tors{(\EC_2)}{p}$ are isomorphic. Then, the sets $\ssprimeset[1]$ and $\ssprimeset[2]$ of primes of supersingular reduction for $\EC_1$ and $\EC_2$ coincide. Given a vector $\multisign\in\{+,-\}^{\ordssprimeset}$, assume that both curves satisfy~\ref{hyp_2} and let $\iwmu_{\EC_j}$ be the Iwasawa $\iwmu[\empty]$-invariants of $\KobdualSel*{\tors{(\EC_j)}{p^\infty}}{\basecyc}$, for $j=1,2$. Then 
\begin{equation}\label{eq:thm:mu=mu:equality_mu}
\iwmu_{\EC_1}=0\Longleftrightarrow \iwmu_{\EC_2}=0.
\end{equation}
\end{theorem}
\begin{proof} 
Observe first that if $\iwmu_{\EC_j}=0$ for one curve $\EC_j$, then~\ref{conj:A} holds for both curves, thanks to Proposition~\ref{prop:conjectures}. Moreover, Proposition~\ref{prop:raynaud}-\ref{item:prop:raynaud_Sss=Sord} shows that the sets $\ssprimeset$ and $\ordprimeset$ consisting of primes of supersingular (\resp ordinary) reduction for $\EC_1$ and $\EC_2$ coincide.

Fix an isomorphism $\tors{(\EC_1)}{p}\cong\tors{(\EC_2)}{p}$ and consider the maps
\[
\plocallambda{\tors{(\EC_j)}{p}}\colon \HH{1}{\GS[\cyc]}{\tors{(\EC_j)}{p}}\longrightarrow \bigoplus_{\genprimeid\in\pprimeset}\modJJ[\pm][\EC_j]{\genprimeid}\kern-2pt=\kern-5pt\bigoplus_{\ordprimeid{}\in\ordprimeset}\bigoplus_{\GenPrimeId\mid\ordprimeid{}}\HH{1}{\basecyc[\GenPrimeId]}{\tors{(\modEC_j)}{p}}\oplus\bigoplus_{i=1}^{\ordssprimeset}\HH[\big]{1}{\localfield[i]{n}}{\tors{(\EC_j)}{p}}/\image \pmKum{\localfield[i]{n}}{p}
\]
defined before Theorem~\ref{thm:main_result}. For all $\GenPrimeId\mid\ordprimeid{}\in\ordprimeset$, the chosen isomorphism induces an isomorphism between the $\HH{1}{\basecyc[\GenPrimeId]}{\tors{(\modEC_j)}{p}}$ (for $j=1,2$) by to Proposition~\ref{prop:raynaud}-\ref{item:prop:raynaud_modE=modE}. Similarly, at every prime in $\ssprimeset$, Proposition~\ref{prop:universal_sequence} gives an isomorphism between the groups
\[
\HH[\big]{1}{\localfield[]{n}}{\tors{(\EC_j)}{p}}/\image \pmKum{\localfield[]{n}}{p},
\]
for $j=1,2$. It follows that $\plocallambda{\tors{(\EC_1)}{p}}$ is surjective if and only if $\plocallambda{\tors{(\EC_2)}{p}}$ is surjective. Theorem~\ref{thm:main_result} now yields~\eqref{eq:thm:mu=mu:equality_mu}.
\end{proof}
In order to prove the next theorem, we need the following proposition, which is the analogue for the Galois representation $\tors{\EC}{p^\infty}$ of the equivalence between~\ref{thm:main_result:surj} and~\ref{thm:main_result:cotorsion} of Theorem~\ref{thm:main_result}, as well as its Corollary~\ref{cor:no_finite_in_signedSel}. To state them, recall the exact sequence
\begin{equation}\tag{\ref{seq:pmcpt_cyc_infty}}\begin{aligned}
0\longrightarrow\KobSel{\tors{\EC}{p^\infty}}{\basecyc}\longrightarrow \HH{1}{\GS[\cyc]}{\tors{\EC}{p^\infty}}\xrightarrow{\pmlambda{\basecyc}{p^\infty}} \bigoplus_{\genprimeid\in \primeset}&\JJ{\genprimeid}
\longrightarrow\pontr{\KobprojSel{\tate{\EC}}{\basecyc}}\longrightarrow\\
\longrightarrow& \HH{2}{\GS[\cyc]}{\tors{\EC}{p^\infty}}\longrightarrow \bigoplus_{\GenPrimeId\mid\genprimeid\in \primeset}\HH{2}{\basecyc[\GenPrimeId]}{\tors{\EC}{p^\infty}}\longrightarrow 0.
\end{aligned}\end{equation}
\begin{proposition}\label{prop:surj_iff_cotors} Under the standing assumption~\ref{hyp_1}, assume further that~\ref{conj:A} holds. Then~\ref{hyp_2} holds if and only if $\pmlambda{\basecyc}{p^\infty}$ is surjective.
\end{proposition}
\begin{proof} By Proposition~\ref{prop:conjectures},~\ref{conj:A} yields that $\HH{2}{\GS[\cyc]}{\tors{\EC}{p^\infty}}=0$, thus $\coker\bigl(\pmlambda{\basecyc}{p^\infty}\bigr)=\pontr{\KobprojSel{\tate{\EC}}{\basecyc}}$. Thanks to~\eqref{seq:pmcpt_cyc_infty}, we need to show that $\pontr{\KobprojSel{\tate{\EC}}{\basecyc}}=0$ if and only if $\KobdualSel{\tors{\EC}{p^\infty}}{\basecyc}$ is $\IwLambda$-torsion. 
The rank computations performed in Proposition~\ref{prop:ranks} yield that
\[
\corank[\IwLambda]\coker\bigl(\pmlambda{\basecyc}{p^\infty}\bigr)=\corank[\IwLambda]\ker\bigl(\pmlambda{\basecyc}{p^\infty}\bigr)=\corank[\IwLambda]\KobSel{\tors{\EC}{p^\infty}}{\basecyc}.
\]
Therefore,~\ref{hyp_2} is equivalent to the statement that $\pontr{\coker\bigl(\pmlambda{\basecyc}{p^\infty}\bigr)}=\KobprojSel{\tate{\EC}}{\basecyc}$ is $\IwLambda$-torsion. By Lemma~\ref{lemma:duality_lim} this can happen if and only if it is trivial, finishing the proof.
\end{proof}
Before deriving Corollary~\ref{cor:no_finite_in_signedSel}, recall the following result due to Greenberg (see~\cite[pp.~104--105]{Gre99}). We give a different, homological proof that works in a broader context. The second author is grateful to Akhil Matthew for helpful discussions in this regard.

\begin{proposition}\label{prop:greenberg_no_finite} Let $0\to M\to N\to M/N\to 0$ be an exact sequence of $\IwLambda$-modules. Suppose that $M$ is free and that $N$ has no non-zero finite $\IwLambda$-submodules. Then $M/N$ has no non-zero finite $\IwLambda$-submodules.
\end{proposition}
\begin{proof} Let us prove that the maximal finite $\IwLambda$-submodule $W\subseteq M/N$ is zero. Being pseudo-null, it satisfies $\operatorname{Ext}^i_{\IwLambda}(W,\IwLambda)=0$, for $i=0,1$, as shown in~\cite[Proposition~5.5.3(ii)]{NeuSchWin08}. Since $M$ is free, this implies $\operatorname{Ext}^1_{\IwLambda}(W,M)=0$. Applying the functor $\Hom_{\IwLambda}(W,-)$ to the exact sequence in the statement therefore gives an exact sequence
\[
0\longrightarrow \Hom_{\IwLambda}(W,M)\longrightarrow \Hom_{\IwLambda}(W,N)\longrightarrow \Hom_{\IwLambda}(W,M/N)\longrightarrow 0.
\]
As $N$ contains no non-zero finite $\IwLambda$-submodules, the term $\Hom_{\IwLambda}(W,N)$ vanishes, and therefore $\Hom_{\IwLambda}(W,M/N)=0$. This implies the proposition, since $W\subseteq M/N$.
\end{proof}
\begin{corollary}\label{cor:no_finite_in_signedSel} Assume~\ref{hyp_1} and~\ref{hyp_2} as well as~\ref{conj:A}. If the Pontryagin dual $\dualSel*{\tors{\EC}{p^\infty}}{\basecyc}$ of the usual Selmer group does not have any non-zero finite $\IwLambda$-submodule, then the same holds for the Pontryagin dual $\KobdualSel*{\tors{\EC}{p^\infty}}{\basecyc}$ of the multi-signed Selmer group.
\end{corollary}
\begin{proof} 
For every $\genprimeid\in\primeset$, define the usual local conditions as
\[
\JJ[\empty]{\genprimeid}=\bigoplus_{\GenPrimeId\mid\genprimeid}\HH{1}{\genlocal[\GenPrimeId]{\cyc}}{\tors{\EC}{p^\infty}}/\image(\fullKum{\genlocal[\GenPrimeId]{\cyc}}{p^\infty}):
\]
the group $\JJ[\empty]{\genprimeid}$ is a quotient of $\JJ{\genprimeid}$, and it actually coincides with it for all $\genprimeid\in\primeset\setminus\ssprimeset$. Consider the commutative triangle
\[
\xymatrix{\HH{1}{\GS[\cyc]}{\tors{\EC}{p^\infty}}\ar@{->}^{\pmlambda{\basecyc}{p^\infty}}[0,2]\ar@{->}_{\fullambda{\basecyc}{p^\infty}}[1,2]&&\displaystyle{\bigoplus_{\genprimeid\in\primeset}\JJ{\genprimeid}}\ar@{->>}[1,0]^\eta\\
&&\displaystyle{\bigoplus_{\genprimeid\in\primeset}\JJ[\empty]{\genprimeid}}
}
\]
where $\eta$ is the canonical surjection. By Proposition~\ref{prop:surj_iff_cotors}, the morphism $\pmlambda{\basecyc}{p^\infty}$ is surjective, and therefore there is a short exact sequence
\begin{equation}\label{seq:pm_full}
0\longrightarrow \ker\pmlambda{\basecyc}{p^\infty}\longrightarrow \ker\fullambda{\basecyc}{p^\infty}\longrightarrow \ker\eta\longrightarrow 0.
\end{equation}
By definition, $\ker\pmlambda{\basecyc}{p^\infty}=\KobSel{\tors{\EC}{p^\infty}}{\basecyc}$ and $\ker\fullambda{\basecyc}{p^\infty}=\fullSel{\tors{\EC}{p^\infty}}{\basecyc}$, while
\[
\ker\eta=\bigoplus_{i=1}^{\ordssprimeset}\image(\fullKum{\ssloccyc[i]}{p^\infty})/\image(\pmKum{\ssloccyc[i]}{p^\infty}).
\]
Observe now that $\image(\fullKum{\ssloccyc[i]}{p^\infty})=\HH{1}{\genlocal[\GenPrimeId]{\cyc}}{\tors{\EC}{p^\infty}}$ by~\cite[Proposition~4.8]{CoaGre96}, so that $\pontr{(\ker\eta)}$ is $\IwLambda$-free by~\cite[Proposition~3.32]{KitOts18}. The corollary follows from Proposition~\ref{prop:greenberg_no_finite} applied to the Pontryagin dual of~\eqref{seq:pm_full}.
\end{proof}
\begin{remark}\label{rmk:criteria_sel}
The above result has been obtained by Kitajima--Otsuki in~\cite[Theorem~4.8]{KitOts18} in the special case when $\multisign=\{+,+,\dots,+\}$ or $\multisign=\{-,-,\dots,-\}$ and under the more restrictive hypothesis that both $\KobdualSel[+]{\tors{\EC}{p^\infty}}{\basecyc}$ \emph{and} $\KobdualSel[-]{\tors{\EC}{p^\infty}}{\basecyc}$ are torsion $\IwLambda$-modules. Our approach through the exact sequence~\eqref{seq:pmcpt_cyc_infty} allows us to make the weaker assumption that just one of these modules be torsion. On the other hand, assuming that both $\KobdualSel[+]{\tors{\EC}{p^\infty}}{\basecyc}$ \emph{and} $\KobdualSel[-]{\tors{\EC}{p^\infty}}{\basecyc}$ are torsion $\IwLambda$-modules, they prove in~\cite[Theorem~4.5]{KitOts18} that the Pontryagin dual  $\dualSel{\tors{\EC}{p^\infty}}{\basecyc}$ of the usual Selmer group does not contain any non-zero, finite, $\IwLambda$-submodule. 

Other criteria for the vanishing of the maximal finite $\IwLambda$-submodule of $\dualSel{\tors{\EC}{p^\infty}}{\basecyc}$ have been given by Lei and the second author in~\cite[Theorem~2.14 and Theorem~4.8]{LeiSuj20}.
\end{remark}

Along the lines of Kim's and Hatley--Lei's results quoted above, when the $\iwmu$-invariants of two residually isomorphic elliptic curves vanish, we can relate their $\iwlambda$-invariants: this is the main content of Theorem~\ref{thm:lambda_inv} below. We denote by
\begin{equation}\label{eq:def_rho}
\seldim{\EC}=\dim_{\finitefield}\dualpmSel{\tors{\EC}{p}}{\basecyc}
\end{equation}
the $\finitefield$-dimension of $\dualpmSel{\tors{\EC}{p}}{\basecyc}$. It clearly depends only upon the isomorphism class of $\tors{\EC}{p}$ and not on $\EC$ itself.
\begin{theorem}\label{thm:lambda_inv}
Let $\EC_1,\EC_2$ be two elliptic curves defined over $\basefield$, satisfying the hypotheses of Theorem~\ref{thm:mu=mu}. Assume that their Iwasawa $\iwmu$-invariants vanish and suppose further that the Pontryagin duals $\dualSel*{\tors{(\EC_j)}{p^\infty}}{\basecyc}$ of the usual Selmer groups do not have any non-zero finite $\IwLambda$-submodule. Then,
\begin{equation}\label{eq:thm:mu=mu:lambda}
\iwlambda_{\EC_j}=\seldim{\empty}+\defect{\EC_j}
\end{equation}
where $\defect{\EC_j}$ is as in Definition~\ref{def:defect} and $\seldim{\empty}:=\seldim{\EC_1}=\seldim{\EC_2}$ is as in~\eqref{eq:def_rho}.
\end{theorem}
\begin{remark}\label{rmk:difference_lambda_+/-} 
The hypothesis that the Pontryagin dual of the signed Selmer groups have no non-zero finite $\IwLambda$-submodule, which we deduce from the analogous result for the usual Selmer group, is known to hold in many cases: see, for instance,~\cite[Theorem~4.5]{KitOts18} and~\cite[Theorem~3.1]{HatLei19}.
\end{remark}
\begin{remark}
As already observed in Corollary~\ref{cor:H[p]=H(p)_infinite_selmer}, the quantity $\defect{\EC}$ is independent of the vector~$\multisign$. In particular, for elliptic curves satisfying~\ref{hyp_1} and~\ref{hyp_2} for two vectors $\multisign_1$ and $\multisign_2$, and such that their multi-signed invariants $\iwmu[\multisign_1]$ and $\iwmu[\multisign_2]$ are both $0$, the difference $\iwlambda[\multisign_1]-\iwlambda[\multisign_2]$ of the multi-signed $\iwlambda[\empty]$-invariants depends only on the isomorphism class of the residual representation.

The concrete manner in which this result will be applied later, is the following. Suppose we are given a family of elliptic curves (satisfying~\ref{hyp_1}) with the property that their residual representations are isomorphic. Given a vector $\multisign$, suppose that all members in the family satisfy~\ref{hyp_2}. A consequence of Theorem~\ref{thm:mu=mu} is that if one member $\genEC$ in the family satisfies $\iwmu[\multisign]_\genEC=0$, then for all other members $\EC$ in the family, we obtain $\iwmu[\multisign]_\EC=0$. Moreover, Theorem~\ref{thm:lambda_inv} shows that if $\genEC$ satisfies $\iwmu[\multisign_1]_\genEC=\iwmu[\multisign_2]_\genEC=0$ for two vectors $\multisign_1,\multisign_2$, the difference of multi-signed Iwasawa invariants
\[
\iwlambda[\multisign_1]_{\EC}-\iwlambda[\multisign_2]_{\EC}=\seldim[\multisign_1]{\genEC}-\seldim[\multisign_2]{\genEC}
\]
is constant. In particular, if $\seldim[\multisign_1]{\genEC}=\seldim[\multisign_2]{\genEC}$, then $\iwlambda[\multisign_1]_{\EC}=\iwlambda[\multisign_2]_{\EC}$ for all curves $\EC$ in the family. This hints at an appropriate generalisation of our results for Coleman families and we hope to investigate it in the future.
\end{remark}
\begin{proof} Corollary~\ref{cor:no_finite_in_signedSel} implies that the Pontryagin duals $\KobdualSel*{\tors{(\EC_j)}{p^\infty}}{\basecyc}$ of the multi-signed Selmer groups do not have any non-zero finite $\IwLambda$-submodule, for $j=1,2$. Further, they are $\IwLambda$-torsion thanks to hypothesis~\ref{hyp_2}. Since $\iwmu_{\EC_1}=\iwmu_{\EC_2}=0$,
\[
\iwlambda_{\EC_j}=\operatorname{length}\Bigl(\KobdualSel{\tors{(\EC_j)}{p^\infty}}{\basecyc}/p\KobdualSel{\tors{(\EC_j)}{p^\infty}}{\basecyc}\Bigr).
\]
On the other hand, taking Pontryagin duals in Corollary~\ref{cor:H[p]=H(p)_infinite_selmer} gives an exact sequence
\begin{equation}\label{seq:W_j}
 V_j\longhookrightarrow\pontr{\Bigl(\tors{\KobSel{\tors{(\EC_j)}{p^\infty}}{\basecyc}}{p}\Bigr)}=\KobdualSel{\tors{(\EC_j)}{p^\infty}}{\basecyc}/p\KobdualSel{\tors{(\EC_j)}{p^\infty}}{\basecyc}\longtwoheadrightarrow\dualpmSel{\tors{(\EC_j)}{p}}{\basecyc}
\end{equation}
where $V_j$ is an $\finitefield$-vector space of finite dimension. Since we are assuming $\iwmu_{\EC_j}=0$, Theorem~\ref{thm:main_result} implies that $\pmlambda{\basecyc}{\tors{(\EC_j)}{p}}$ is surjective and therefore, again by Corollary~\ref{cor:H[p]=H(p)_infinite_selmer}, we have $\dim_{\finitefield}V_j=\defect{\EC_j}$. Taking lengths in~\eqref{seq:W_j} gives
\[
\iwlambda_{\EC_j}=\seldim{\empty}+\defect{\EC_j}.
\]
and this finishes the proof.
\end{proof}
In the next two corollaries, we consider the main setting of Theorem~\ref{thm:mu=mu}. Thus, let $\EC_1,\EC_2$ be two elliptic curves defined over $\basefield$ satisfying hypotheses~\ref{hyp_1} and~\ref{hyp_2}, and such that the residual Galois representations $\tors{(\EC_1)}{p}$ and $\tors{(\EC_2)}{p}$ are isomorphic, so $\seldim{\EC_1}=\seldim{\EC_2}:=\seldim{\empty}$. By Proposition~\ref{prop:raynaud}-\ref{item:prop:raynaud_Sss=Sord}, the set $\ordprimeset$ of $p$-adic primes where the curves have good, ordinary reduction, coincide. Further, we assume that $\iwmu_{\EC_1}=0$, which is equivalent to assuming $\iwmu_{\EC_2}=0$ by Theorem~\ref{thm:mu=mu}. Moreover, this implies that~\ref{conj:A} holds for both curves, again by Theorem~\ref{thm:mu=mu}. Finally, suppose that the Pontryagin duals of the usual Selmer groups do not contain any non-trivial, finite, $\IwLambda$-submodule., so that Theorem~\ref{thm:lambda_inv} applies.
\begin{corollary} Let $\badprimeset[1]$ and $\badprimeset[2]$, be the sets of primes of bad reduction for $\EC_1$ and $\EC_2$, respectively. 
If, for both indices $j\in\{1,2\}$, we have $\tors{\EC_j(\genlocal{\genprimeid})}{p}=0$ for all $\genprimeid\in\ordprimeset\cup\badprimeset[j]$, then
\[
\iwlambda_{\EC_1}=\iwlambda_{\EC_2}=\seldim{\empty}.
\]
\end{corollary}
\begin{proof} Recall from Definition~\ref{def:defect} that
\begin{equation}\label{eq:defect_j}
\defect{\EC_j}=\sum_{\badprimeid{\empty}\in\badprimeset[j]}\ordcycprimes{\badprimeid{\empty}}\cdot\dim_{\finitefield}\tors{\EC_j(\cyclocal{\BadPrimeId{\empty}})}{p}+\sum_{\ordprimeid{\empty}\in\ordprimeset}\ordcycprimes{\ordprimeid{\empty}}\dim_{\finitefield}\tors{\modEC_j(\finitefield[\ordprimeid{\empty}])}{p}
\end{equation}
where $\BadPrimeId{\empty}$ is a prime in $\basecyc$ above $\badprimeid{\empty}$.

The same argument as in the proof of Corollary~\ref{cor:H[p]=H(p)_infinite_selmer} shows that the condition $\tors{\EC_j(\genlocal{\badprimeid{\empty}})}{p}=0$ is equivalent to $\tors{\EC_j(\genlocal[\cyc]{\BadPrimeId{\empty}})}{p}=0$ for all $\BadPrimeId{\empty}\mid\badprimeid{\empty}$. In particular, the hypothesis of the corollary imply
$\tors{\EC_j(\cyclocal{\BadPrimeId{\empty}})}{p}$ for all $\BadPrimeId{\empty}\mid\badprimeid{\empty}$. Similarly, there are surjections
\[
\tors{\EC_j(\genlocal{\ordprimeid{\empty}})}{p}\twoheadrightarrow \tors{\modEC_j(\finitefield[\ordprimeid{\empty}])}{p}
\]
so $\tors{\EC_j(\genlocal{\ordprimeid{\empty}})}{p}=0$ implies $\tors{\modEC_j(\finitefield[\ordprimeid{\empty}])}{p}=0$.

Hence all terms in~\eqref{eq:defect_j} vanish and $\defect{\EC_1}=\defect{\EC_2}=0$. The corollary follows from~\eqref{eq:thm:mu=mu:lambda}.
\end{proof}
Recall that, given a prime $\genprimeid\in\primeset$, we denote by $\ordcycprimes{\genprimeid}$ the number of primes $\GenPrimeId\mid\genprimeid$ in $\basecyc$.
\begin{corollary}\label{cor:CM} Suppose that $\EC_1$ is a CM curve. Then
\begin{align*}
\iwlambda_{\EC_2}&=\left(\seldim{\empty}+\sum_{\ordprimeid{\empty}\in\ordprimeset}\ordcycprimes{\ordprimeid{\empty}}\dim_{\finitefield}\tors{\modEC_1(\finitefield[\ordprimeid{\empty}])}{p}\right)+\sum_{\phantom{\kern.3em}\badprimeid{}\in\badprimeset[2]}\ordcycprimes{\badprimeid{}}\cdot \dim_{\finitefield}\tors{\EC_2(\genlocal{\badprimeid{}})}{p}
\end{align*}
\end{corollary}
\begin{remark} The interest of Corollary~\ref{cor:CM} lies in the fact that the quantity in parenthesis is constant along families with isomorphic residual representation at $p$. Moreover, the final sum in the right-hand side only depends on the groups $\tors{\EC_2(\genlocal{\badprimeid{}})}{p}$ (for $\badprimeid{}\in\badprimeset[2]$) and not on the behaviour of $p$-torsion along the local cyclotomic $\ZZ_p$-towers. As we shall see in the proof, the corollary still holds only assuming that the image of $\Gal{\overline{\basefield}}{\basefield}$ inside $\operatorname{Aut}\bigl(\tors{(\EC_1)}{p}\bigr)\subseteq\operatorname{GL}_2(\finitefield)$ is contained in the normalizer of a Cartan subroup, which is certainly the case when $\EC_1$ is CM.
\end{remark}
\begin{proof}
Since $\EC_1$ is CM and $p\geq 3$, the image of $\Gal{\overline{\basefield}}{\basefield}$ inside $\operatorname{Aut}\bigl(\tors{(\EC_1)}{p}\bigr)\subseteq\operatorname{GL}_2(\finitefield)$ is contained in the normalizer of a Cartan subroup. In particular, it contains no element of order $p$, and the same holds for the image of $\Gal{\overline{\basefield}}{\basefield}$ inside $\operatorname{Aut}\bigl(\tors{(\EC_2)}{p}\bigr)$ because the representations are isomorphic. It follows that, for all $\BadPrimeId{}\mid\badprimeid{}$, the pro-$p$-group $\Gamma=\Gal{\genlocal[\cyc]{\BadPrimeId{}}}{\genlocal{\badprimeid{}}}$ acts trivially on $\tors{\EC_2(\genlocal[\cyc]{\BadPrimeId{}})}{p}$, and $\dim_{\finitefield}\tors{\EC_2(\genlocal{\badprimeid{}})}{p}=\dim_{\finitefield}\tors{\EC_2(\cyclocal{\BadPrimeId{}})}{p}$. The corollary follows from~\eqref{eq:thm:mu=mu:lambda}, combined with Definition~\ref{def:defect}.
\end{proof}
\section{Numerical Examples}\label{section:numerical_examples}
Our class of examples comes from the work~\cite{RubSil95}. Both for $p=3$ and $p=5$, Rubin and Silverberg define, for each $D\not\equiv 0\pmod{p}$, a family parametrised by\footnote{Actually, the parameters in the families can vary in $\QQ$, but are required to be $p$-integral to define curves with good reduction at $p$. In our examples, we will restrict to $t\in\ZZ$} $t\in\ZZ$. All curves in the families have good, supersingular reduction at $p$ and isomorphic residual Galois representations. In particular, the reduction type is constant along families and, since all curves are defined over $\QQ$, in all cases $\ordprimeset=\emptyset$. Finally, observe that Rubin--Silverberg's construction shows that all family contain a CM member, and so Corollary~\ref{cor:CM} applies. Since the field of definition of all curves is $\QQ$, we have $\ordssprimeset=1$, allowing for simplified notation; hence, write $\pm$ for the generic vector $\multisign\in\{+,-\}$. For all choices of $(p,D)$, the strategy will be to
\begin{enumerate}[label=\arabic*.]
\item Find one curve $\genEC$ in the family for which the Iwasawa invariants $\iwmu[\pm]_{\genEC}$ and $\iwlambda[\pm]_{\genEC}$ have been computed in~\cite{LMFDB} and such that $\iwmu[+]_{\genEC}=\iwmu[-]_{\genEC}=0$. In practice, we take for $\genEC$ the CM curve corresponding to the parameter $t=0$.
\item Apply Theorem~\ref{thm:mu=mu} (see in particular Remark~\ref{rmk:difference_lambda_+/-}) to deduce that $\iwmu[\pm]=0$ for all other members in the family. In particular,~\ref{conj:A} holds for the whole family, by Proposition~\ref{prop:conjectures}.
\item Deduce from Theorem~\ref{thm:lambda_inv} (which can be applied thanks to \cite[Theorem~4.5]{KitOts18}) that $\seldim[\pm]{\genEC}=\iwlambda[\pm]_{\genEC}-\defect{\genEC}$ for $\genEC$, and set $\seldim[\pm]{\empty}:=\seldim[\pm]{\genEC}$.
\item By Corollary~\ref{cor:CM}, we obtain
\[
\iwlambda[\pm]_{\EC}=\seldim[\pm]{\empty}+\defect{\EC}=\seldim[\pm]{\empty}+\sum_{\ell\in\badprimeset}\ordcycprimes{\ell}\dim_{\finitefield}\tors{\EC(\QQ_\ell)}{p}
\]
for all $\EC$ in the family.
\item\label{step:hasse} The key step is to find elliptic curves $\EC$ in the family satisfying $a_p(\EC)=0$, to ensure that~\ref{hyp_1} holds. Note that this is only needed when $p=3$, because when $p=5$ the condition $a_5(\EC)=0$ is automatically satisfied by the Hasse bound. Since all our examples are defined over $\QQ$,~\ref{hyp_2} is always satisfied by~\cite[Theorem~1.2]{Kob03} and the usual Selmer group never contain non-trivial, finite, $\IwLambda$-submodules by~\cite[Theorem~4.5]{KitOts18} combined with Corollary~\ref{cor:no_finite_in_signedSel}.
\item Choosing any curve as in~\eqref{step:hasse}, we compute the $\finitefield$-dimension of $\tors{\EC(\QQ_\ell)}{p}$ at all primes $\ell\in\badprimeset$, together with the number of primes in $\Qcyc$ above $\ell$, to find the numerical value of $\defect{\EC}$ and hence of $\iwlambda[\pm]_{\EC}$.
\end{enumerate}
We will consider the families attached to $D=1,-1$ for $p=3$, and the families attached to $D=3$ and $D=14$ for $p=5$. Our source of numerical data is~\cite{LMFDB}. Labels of elliptic curves follow Cremona's tables as in~\cite{LMFDB}, when available (\ie for discriminant less than $500.000$ as per October~$2019$). The computations have been made in \texttt{SAGE}\footnote{
We used commands \texttt{E.q\_expansion(4)} to compute $a_3$ and \texttt{E(0).division\_points(p)} to compute torsion points.
}.
\subsection{\texorpdfstring{$\boldsymbol{p=3}$}{p=3}}
\subsubsection{$D=1$} Setting $t=0$ we obtain the CM curve $\genEC=32a2$ given by $y^2=x^3-x$. It satisfies $\iwmu[\pm]_\genEC=\iwlambda[\pm]_\genEC=0$ and this is in accordance with the fact that we found $\defect{\genEC}=0$: indeed, $\badprimeset[\genEC]=\{2\}$ and $\tors{\genEC(\QQ_2)}{3}=0$. Moreover, $a_3(\genEC)=0$, and we obtain $\seldim[\pm]{\empty}=0$. The curves corresponding to $t=1$ and $t=2$ are, respectively, $\EC_1=352f1$ and $\EC_2=16096h1$: since $a_3(\EC_1)=-3$ and $a_3(\EC_2)=3$, we discard them.

The curve corresponding to $t=3$ is $\EC_3=18784b1$, and $a_3(\EC_3)=0$. Its Iwasawa invariants are available on~\cite{LMFDB}, and indeed $\iwmu[\pm]_{\EC_3}=0$. The primes of bad reduction are $\badprimeset=\{2,587\}$. We found $\tors{\EC_3(\QQ_2)}{3}=0$ and $\tors{\EC_3(\QQ_{587})}{3}=\ZZ/3\ZZ$; since $587$ is a generator of $\ZZ/9\ZZ$, it is totally inert in $\Qcyc/\QQ$, so $\ordcycprimes{357}=1$. Formula~\eqref{eq:thm:mu=mu:lambda} gives $\iwlambda[\pm]_{\EC_3}=1$, in accordance with the numerical value found in~\cite{LMFDB}.

To show a somehow extreme example, consider $t=18$. It satisfies $a_3(\EC_{18})=0$ and its conductor is $90.885.856=2^5\cdot 2840183$. The dimensions of its local $3$-torsion are
\[
\dim_{\finitefield[3]}\Bigl(\tors{\EC_{18}(\QQ_\ell)}{3}\Bigr)=\begin{cases}
0&\text{ for }\ell=2\\
1&\text{ for }\ell=2840183\end{cases}
\]
The multiplicative order of $2840183$ modulo $3^7$ being $6$, we deduce $\ordcycprimes{2840183}=3^6$, whence $\iwlambda[\pm]_{\EC_{18}}=729$, and $\iwmu[\pm]_{\EC_{18}}=0$ by Theorem~\ref{thm:mu=mu}. It is relevant here to note Kim's observation that under these assumptions the Iwasawa $\iwlambda[\pm]$-invariants can be arbitrarily large in the family (see~\cite[p.~190]{Kim09}), although he does not produce explicit examples. Note also that these Iwasawa invariants are not available on~\cite{LMFDB}.
\subsubsection{$D=-1$} In this case, the CM curve for $t=0$ is $\genEC=64a4$ given by $y^2=x^3+x$. Again, $\iwmu[\pm]_{\genEC}=\iwlambda[\pm]_{\genEC}=0=a_3(\genEC)$. We computed the defect and found $\defect{\genEC}=0$, since $\badprimeset[\genEC]=\{2\}$ and $\tors{\genEC(\QQ_2)}{3}=0$. We obtain $\seldim[\pm]{\empty}=0$. The curves corresponding to the parameters $t=2,4,5$ are, respectively, $\EC_2=22976p1,\EC_4=423872t1$ and $\EC_{5}=131392f1$. They all exist in~\cite{LMFDB}, and have $a_3(\EC_i)=0$, but the Iwasawa invariants are available only for $\EC_2$ and $\EC_5$: they read $\iwlambda[\pm]_{\EC_2}=3,\iwlambda[\pm]_{\EC_5}=0$. This is in accordance with formula~\eqref{eq:thm:mu=mu:lambda}: indeed, $\badprimeset[\EC_2]=\{2,359\},\badprimeset[\EC_5]=\{2,2053\}$ and
\[
\dim_{\finitefield[3]}\Bigl(\tors{\EC_{2}(\QQ_\ell)}{3}\Bigr)=\begin{cases}
0&\text{ for }\ell=2\\
1&\text{ for }\ell=359\\
\end{cases}\qquad
\dim_{\finitefield[3]}\Bigl(\tors{\EC_{5}(\QQ_\ell)}{3}\Bigr)=\begin{cases}
0&\text{ for }\ell=2\\
0&\text{ for }\ell=2053\\
\end{cases}
\]
This immediately implies $\defect{\EC_5}=0$, so $\iwlambda[\pm]_{\EC_5}=0$. As $359$ has order $6$ modulo ${27}$, we obtain $\ordcycprimes{359}=3$, whence $\defect{\EC_2}=\iwlambda[\pm]_{\EC_5}=3$. The curve $\EC_4$ can be treated analogously, since $\badprimeset[\EC_4]=\{2,37,179\}$ and
\[
\dim_{\finitefield[3]}\Bigl(\tors{\EC_{4}(\QQ_\ell)}{3}\Bigr)=\begin{cases}
0&\text{ for }\ell=2\\
2&\text{ for }\ell=37\\
1&\text{ for }\ell=179\\
\end{cases}
\]
Further, $\ordcycprimes{37}=\ordcycprimes{179}=3$, whence $\iwlambda[\pm]_{\EC_4}=9$, a value which is not available on~\cite{LMFDB}. Also, all curves satisfy $\iwmu[\pm]=0$ by Theorem~\ref{thm:mu=mu}.

We finish this series of examples with the curve $\EC_{149}$ for $t=149$. Its conductor is $106.459.833.664=2\cdot 1663434901$, so it has no label in Cremona's tables, but we can compute $a_3(\EC_{149})=0$. We found $\tors{\EC_{149}(\QQ_2)}{3}=\tors{\EC_{149}(\QQ_{1663434901})}{3}=0$, whence $\iwmu[\pm]_{\EC_{149}}=\iwlambda[\pm]_{\EC_{149}}=0$.
\subsection{\texorpdfstring{$\boldsymbol{p=5}$}{p=5}}
\subsubsection{$D=3$}
The CM curve corresponding to $t=0$ is $\genEC=3888s1$, given by $y^2=x^3+48$. Its Iwasawa invariants are computed in~\cite{LMFDB} and $\iwmu[\pm]_{\genEC}=0,\iwlambda[\pm]_{\genEC}=1$. To find $\seldim[\pm]{\empty}=\seldim[\pm]{\genEC}$, we need to compute $\defect{\genEC}$. The primes of bad reduction are $\badprimeset[\genEC]=\{2,3\}$ and $\tors{\genEC(\QQ_\ell)}{3}=0$ for both $\ell\in\badprimeset[\genEC]$, so $\defect{\genEC}=0$ and $\seldim[\pm]{\empty}=1$. The conductors of $\EC_t$ for $t\in[-5,15]$ have orders of magnitude between $10^7$ and $10^{20}$ (except for $\EC_0=\genEC$), so these curves are not implemented in~\cite{LMFDB}. Computing Iwasawa invariants through formula~\eqref{eq:thm:mu=mu:lambda} is almost immediate. As an example, we compute them for the curves $\EC_6$ and $\EC_{14}$ corresponding to $t=6$ and $t=14$, respectively. First, we immediately obtain from Theorem~\ref{thm:mu=mu} that $\iwmu[\pm]_{\EC_6}=\iwmu[\pm]_{\EC_{14}}=0$.

The conductor of $\EC_6$ is $16.847.046.490.346.928=2^4\cdot 3^5\cdot 4333088089081$. The curve has no $\QQ_\ell$-rational $5$-torsion points for any of the primes $\ell\in\{2,3,4333088089081\}$, so $\defect{\EC_6}=0$. It follows that $\iwlambda[\pm]_{\EC_6}=\seldim[\pm]{\empty}=1$. The conductor of $\EC_{14}$ is $445.766.016.078.830.163.888=2^4\cdot 3^5 \cdot 29\cdot 602279\cdot 6564248011$ and $\EC_{14}$ does not have neither $\QQ_2$-rational nor $\QQ_3$-rational $5$-torsion points. On the other hand,
\[
\dim_{\finitefield[5]}\Bigl(\tors{\EC_{14}(\QQ_\ell)}{5}\Bigr)=\begin{cases}
1&\text{ for }\ell=29\\
1&\text{ for }\ell=602279\\
2&\text{ for }\ell=6564248011\\
\end{cases}
\]
Further, computing multiplicative orders modulo $25$, we find $\ordcycprimes{\ell}=1$ for all $\ell\in\{29, 602279,6564248011\}$. It follows that $\defect{\EC_{14}}=4$ and $\iwlambda[\pm]_{\EC_{14}}=\defect{\EC_{14}}+\seldim[\pm]{\empty}=5$.
\subsubsection{$D=14$}
We finish with an example where $\iwlambda[+]\neq\iwlambda[-]$. Take $D=14$, so that the CM member for $t=0$ is $\genEC=28224dj1$, given by $y^2=x^2+224$. Its Iwasawa invariants are $\iwmu[\pm]_{\genEC}=0$, and $\iwlambda[+]_{\genEC}=3,\iwlambda[-]_{\genEC}=1$. The conductor of $\genEC$ is $28224=2^{6}\cdot 3^{2} \cdot 7^{2}$ and we compute as above that $\defect{\genEC}=0$, so $\seldim[+]{\genEC}=\seldim[+]{\empty}=3,\seldim[-]{\genEC}=\seldim[-]{\empty}=1$. As observed in Remark~\ref{rmk:difference_lambda_+/-}, all members $\EC_t$ in this family satisfy $\iwlambda[+]_{\EC}-\iwlambda[-]_{\EC}=2$, together with $\iwmu[\pm]=0$. Again, the conductors grow very fast with $t$ and we could not find any curve in the family for which data are available on~\cite{LMFDB}. As examples, we consider the curves for $t=6$ and $t=8$.
The first has $\badprimeset[\EC_6]=\{2,3,7,22621,92081500261\}$ and there are no $\QQ_\ell$-rational $5$-torsion points at $\ell\in\badprimeset[\EC_6]$ except for $\ell=92081500261$, where $\dim_{\finitefield[5]}\bigl(\tors{\EC_{6}(\QQ_\ell)}{5}\bigr)=2$. Since $\ordcycprimes{92081500261}=1$, we find $\defect{\EC_6}=2$ and 
\[
\iwlambda[+]_{\EC_6}=5\qquad\text{ and }\qquad \iwlambda[-]_{\EC_6}=3.
\]
Finally, we consider the curve for $t=8$, which has no $\QQ_\ell$-rational $5$-torsion point at any of the primes $\ell\in\badprimeset[\EC_8]=\{2,3,7,10861,642211,9447511\}$. It follows that $\defect{\EC_8}=0$ and 
\[
\iwlambda[+]_{\EC_6}=3\qquad\text{ and }\qquad \iwlambda[-]_{\EC_6}=1.
\]

\bibliographystyle{amsalpha}
\bibliography{arXiv_bib}

\providecommand{\bysame}{\leavevmode\hbox to3em{\hrulefill}\thinspace}
\providecommand{\MR}{\relax\ifhmode\unskip\space\fi MR }
\providecommand{\MRhref}[2]{%
  \href{http://www.ams.org/mathscinet-getitem?mr=#1}{#2}
}
\providecommand{\href}[2]{#2}
\begin{thebibliography}{{LMF}13}

\bibitem[CG96]{CoaGre96}
J.~Coates and R.~Greenberg, \emph{Kummer theory for abelian varieties over
  local fields}, Invent. Math. \textbf{124} (1996), no.~1-3, 129--174.
  \MR{1369413}

\bibitem[CS05]{CoaSuj05}
John Coates and Ramdorai Sujatha, \emph{Fine {S}elmer groups of elliptic curves
  over $p$-adic {L}ie extensions}, Math. Ann. \textbf{331} (2005), no.~4,
  809--839.

\bibitem[CS10]{CoaSuj10}
John Coates and Ramdorai Sujatha, \emph{Galois cohomology of elliptic curves},
  second ed., Published by Narosa Publishing House, New Delhi; for the Tata
  Institute of Fundamental Research, Mumbai, 2010. \MR{3060733}

\bibitem[EPW06]{EmePolWes06}
Matthew Emerton, Robert Pollack, and Tom Weston, \emph{Variation of iwasawa
  invariants in hida families}, Invent. Math. \textbf{163} (2006), no.~3,
  523--580.

\bibitem[Gre89]{Gre89}
Ralph Greenberg, \emph{Iwasawa theory for {$p$}-adic representations},
  Algebraic number theory, Adv. Stud. Pure Math., vol.~17, Academic Press,
  Boston, MA, 1989, pp.~97--137. \MR{1097613}

\bibitem[Gre99]{Gre99}
\bysame, \emph{Iwasawa theory for elliptic curves}, Arithmetic theory of
  elliptic curves ({C}etraro, 1997), Lecture Notes in Math., vol. 1716,
  Springer, Berlin, 1999, pp.~51--144. \MR{1754686}

\bibitem[Gre11]{Gre11}
\bysame, \emph{Iwasawa theory, projective modules, and modular
  representations}, Mem. Amer. Math. Soc. \textbf{211} (2011), no.~992, vi+185.
  \MR{2807791}

\bibitem[GV00]{GreVat00}
Ralph Greenberg and Vinayak Vatsal, \emph{On the iwasawa invariants of elliptic
  curves}, Invent. Math. \textbf{142} (2000), no.~1, 17--63.

\bibitem[Hac11]{Hac11}
Yoshitaka Hachimori, \emph{Iwasawa {$\lambda$}-invariants and congruence of
  {G}alois representations}, J. Ramanujan Math. Soc. \textbf{26} (2011), no.~2,
  203--217. \MR{2816789}

\bibitem[HL19]{HatLei19}
Jeffrey Hatley and Antonio Lei, \emph{Arithmetic properties of signed selmer
  groups at non-ordinary primes}, Annales de l'Institut Fourier \textbf{69}
  (2019), no.~3, 1259--1294 (en).

\bibitem[IP06]{IovPol06}
Adrian Iovita and Robert Pollack, \emph{Iwasawa theory of elliptic curves at
  supersingular primes over {$\mathbb{Z}_p$}-extensions of number fields}, J.
  Reine Angew. Math. \textbf{598} (2006), 71--103. \MR{2270567}

\bibitem[Jan14]{Jan14}
Uwe Jannsen, \emph{A spectral sequence for {I}wasawa adjoints}, M\"unster J.
  Math. \textbf{7} (2014), no.~1, 135--148. \MR{3271243}

\bibitem[Kat04]{Kat04}
Kazuya Kato, \emph{$p$-adic {H}odge theory and values of zeta functions of
  modular forms}, Asterisque (2004), no.~295, ix, 117--290.

\bibitem[Kid18]{Kid18}
Keenan Kidwell, \emph{On the structure of {S}elmer groups of {$p$}-ordinary
  modular forms over {${\bf{Z}}_p$}-extensions}, J. Number Theory \textbf{187}
  (2018), 296--331. \MR{3766913}

\bibitem[Kim09]{Kim09}
Byoung~Du Kim, \emph{The {I}wasawa invariants of the plus/minus {S}elmer
  groups}, Asian J. Math. \textbf{13} (2009), no.~2, 181--190. \MR{2559107}

\bibitem[Kim13]{Kim13}
\bysame, \emph{The plus/minus {S}elmer groups for supersingular primes}, J.
  Aust. Math. Soc. \textbf{95} (2013), no.~2, 189--200. \MR{3142355}

\bibitem[Kim18]{Kim18}
\bysame, \emph{Ranks of the rational points of abelian varieties over ramified
  fields, and {I}wasawa theory for primes with non-ordinary reduction}, J.
  Number Theory \textbf{183} (2018), 352--387. \MR{3715241}

\bibitem[KO18]{KitOts18}
Takahiro Kitajima and Rei Otsuki, \emph{On the plus and the minus selmer groups
  for elliptic curves at supersingular primes}, Tokyo J. Math. \textbf{41}
  (2018), no.~1, 273--303.

\bibitem[Kob03]{Kob03}
Shin{-}ichi Kobayashi, \emph{Iwasawa theory for elliptic curves at
  supersingular primes}, Invent. Math. \textbf{152} (2003), no.~1, 1--36.

\bibitem[LLZ10]{LeiLoeZer10}
Antonio Lei, David Loeffler, and Sarah~Livia Zerbes, \emph{Wach modules and
  {I}wasawa theory for modular forms}, Asian J. Math. \textbf{14} (2010),
  no.~4, 475--528. \MR{2774276}

\bibitem[{LMF}13]{LMFDB}
The {LMFDB Collaboration}, \emph{The {L}-functions and modular forms database},
  \url{http://www.lmfdb.org}, 2013, [Online; accessed October $31^\text{st}$
  2019].

\bibitem[LS18a]{LimSuj18a}
Meng~Fai Lim and Ramdorai Sujatha, \emph{Fine {S}elmer groups of congruent
  {G}alois representations}, J. Number Theory \textbf{187} (2018), 66--91.

\bibitem[LS18b]{LimSuj18b}
\bysame, \emph{On the structure of fine {S}elmer groups and {S}elmer groups of
  {CM} elliptic curves}, in preparation, 2018.

\bibitem[LS20]{LeiSuj20}
Antonio Lei and Ramdorai Sujatha, \emph{On {S}elmer groups in the supersingular
  reduction case}, Tokyo J. Math. (2020), to appear.

\bibitem[Maz72]{Maz72}
Barry Mazur, \emph{Rational points of abelian varieties with values in towers
  of number fields}, Invent. Math. \textbf{18} (1972), 183--266.

\bibitem[NSW08]{NeuSchWin08}
J{\"u}rgen Neukirch, Alexander Schmidt, and Kay Wingberg, \emph{{Cohomology of
  number fields}}, second ed., {Grundlehren der Mathematischen Wissenschaften
  [Fundamental Principles of Mathematical Sciences]}, vol. 323,
  Springer-Verlag, New York, 2008.

\bibitem[Pol03]{Pol03}
Robert Pollack, \emph{On the p-adic l-function of a modular form at a
  supersingular prime}, Duke Math. J. \textbf{118} (2003), no.~3, 523--558.

\bibitem[PR90]{Per90}
Bernadette Perrin-Riou, \emph{Théorie d'{I}wasawa $p$-adique locale et
  globale}, Invent. Math. \textbf{99} (1990), no.~2, 247--292.

\bibitem[PR92]{Per92}
\bysame, \emph{Théorie d'iwasawa et hauteurs p-adiques}, Invent. Math.
  \textbf{109} (1992), no.~1, 137--185.

\bibitem[PR95]{Per95}
Bernadette Perrin-Riou, \emph{Fonctions {$L$} {$p$}-adiques des
  repr\'{e}sentations {$p$}-adiques}, Ast\'{e}risque (1995), no.~229, 198.
  \MR{1327803}

\bibitem[Ray74]{Ray74}
Michel Raynaud, \emph{Sch\'{e}mas en groupes de type {$(p,\dots, p)$}}, Bull.
  Soc. Math. France \textbf{102} (1974), 241--280. \MR{419467}

\bibitem[RS95]{RubSil95}
K.~Rubin and A.~Silverberg, \emph{Families of elliptic curves with constant mod
  {$p$} representations}, Elliptic curves, modular forms, \& {F}ermat's last
  theorem ({H}ong {K}ong, 1993), Ser. Number Theory, I, Int. Press, Cambridge,
  MA, 1995, pp.~148--161. \MR{1363500}

\bibitem[Sch85]{Sch85}
Peter Schneider, \emph{{\$}p{\$}-adic height pairings. ii}, Invent. Math.
  \textbf{79} (1985), no.~2, 329--374.

\bibitem[Ser94]{Ser94}
Jean-Pierre Serre, \emph{Cohomologie galoisienne}, fifth ed., Lecture Notes in
  Mathematics, vol.~5, Springer-Verlag, Berlin, 1994. \MR{1324577}

\bibitem[Sil94]{Sil94}
Joseph~H. Silverman, \emph{Advanced topics in the arithmetic of elliptic
  curves}, Graduate Texts in Mathematics, vol. 151, Springer-Verlag, New York,
  1994. \MR{1312368}

\bibitem[Spr12]{Spr12}
Florian E.~Ito Sprung, \emph{Iwasawa theory for elliptic curves at
  supersingular primes: a pair of main conjectures}, J. Number Theory
  \textbf{132} (2012), no.~7, 1483--1506. \MR{2903167}

\bibitem[Suj10]{Suj10}
R.~Sujatha, \emph{Elliptic curves and {I}wasawa's {$\mu=0$} conjecture},
  Quadratic forms, linear algebraic groups, and cohomology, Dev. Math.,
  vol.~18, Springer, New York, 2010, pp.~125--135. \MR{2648723}

\end{thebibliography}
\end{document}